\documentclass[smallextended]{svjour3}       
\smartqed  
\usepackage{xspace}
\usepackage{algorithm}
\usepackage{subfigure}
\usepackage{amsfonts}
\usepackage{amsmath}
\usepackage{enumerate}
\usepackage{framed}
\usepackage{graphicx}
\usepackage{lscape}
\usepackage{bm}
\usepackage{color}
\usepackage{multicol}

\newcommand{\Dt}{\Delta t}

\newcommand{\figref}[1]{{Figure~\ref{#1}}}

\newcommand{\secref}[1]{{Section~\ref{#1}}}
\renewcommand{\eqref}[1]{{(\ref{#1})}}
 \newcommand{\HH}{\mathbb{H}}
 \newcommand{\N}{\mathbb{N}}
\newtheorem{assumption}[theorem]{Assumption}
\newcommand{\thmref}[1]{{Theorem~\ref{#1}}}
 \newcommand{\lemref}[1]{{Lemma~\ref{#1}}}
 \newcommand{\assref}[1]{{Assumption~\ref{#1}}}
\newcommand{\propref}[1]{{Proposition~\ref{#1}}}
\newcommand{\rmref}[1]{{Remark ~\ref{#1}}}
\newcommand{\Shdt}{S_{h,\Dt}}

\usepackage{graphicx}
\def\makeheadbox{\relax}

\begin{document}

\title{A modified semi--implicit Euler-Maruyama scheme for finite
 element discretization of SPDEs with additive noise
}

\titlerunning{Modified Scheme for Finite Element
Discretization}        

\author{Gabriel J. Lord        \and
        Antoine Tambue 
}


\institute{ G. J. Lord \at
              \\
            Department of Mathematics and the Maxwell Institute for Mathematical Sciences, Heriot Watt University, Edinburgh EH14 4AS, U.K \\
            Tel.: +44-131-451-8196\\
            Fax: +44-131-451-3249\\
          \email{g.j.lord@hw.ac.uk}           
           \and          
           Antoine Tambue (Corresponding author) \at
            The African Institute for Mathematical Sciences(AIMS) of South Africa and Stellenbosh University,\\
            Center for Research in Computational and Applied Mechanics (CERECAM), and Department of Mathematics and Applied Mathematics, University of Cape Town, 7701 Rondebosch, South Africa.\\
             Tel.: +27-785580321\\
   \email{antonio@aims.ac.za, tambuea@gmail.com}           
}

\date{Received: date / Accepted: date}

\maketitle

\begin{abstract}
  We consider the numerical approximation of a general second order
  semi--linear parabolic stochastic partial differential equation
  (SPDE) driven by additive space-time noise. We introduce a new
  modified scheme using a linear functional of the noise with a
  semi--implicit Euler--Maruyama method in time and in space we analyse a
  finite element method (although extension to finite differences or
  finite volumes would be possible). 
  We prove convergence in the root mean square $L^{2}$ norm for
  a diffusion reaction equation and diffusion advection reaction equation.
  We present numerical results for a linear reaction diffusion equation
  in two dimensions as well as a nonlinear example of two-dimensional
 stochastic advection diffusion reaction equation. 
 We see from both the analysis and numerics that the proposed scheme has
 better convergence properties than the standard semi--implicit Euler--Maruyama method.
 
  
\end{abstract}

\section{Introduction}
\label{intro}
We analyse the strong numerical approximation of Ito
stochastic partial differential equations defined in $\Omega\subset\mathbb{R}^{d}$.
 Boundary conditions on the domain $\Omega$
are typically Neumann, Dirichlet or some mixed conditions.
We consider equations of the form
\begin{eqnarray}
\label{adr}
dX=(AX +F(X))dt + d W, \qquad  X(0)=X_{0},\qquad t \in [0, T],\quad T>0
\end{eqnarray}
in a Hilbert space $H= L^{2}(\Omega)$. Here 
$A$ is the generator of an analytic semigroup $S(t):=e^{t A}, t\geq 0$
with eigenfunctions $e_i$ and eigenvalues $\lambda_i$, $i\in\mathbb{N}^{d}$. $F$ 
is a nonlinear function of $X$ and possibly
$\nabla X$. The noise term, $W(x,t)$, is a $Q$-Wiener process that is
white in time and defined on a filtered probability space
$(\mathbb{D},\mathcal{F},\mathbb{P},\left\lbrace  F_{t}\right\rbrace_{t\geq 0})$.
 We assume that the noise can be represented as 
 \begin{eqnarray}
  \label{eq:W1}
  W(x,t)=\underset{i \in  \mathbb{N}^{d}}{\sum}\sqrt{q_{i}}e_{i}(x)\beta_{i}(t), 
\end{eqnarray}
where  $q_i$, $i\in \mathbb{N}^{d}$ are respectively the eigenvalues and the eignfunctions  of $Q$,
 and $\beta_{i}$ are
independent and identically distributed standard Brownian motions.
Precise assumptions on $A$, $F$ and $W$ are given in
\secref{sec:NumericalScheme} and, under these type of technical assumptions,
it is well known (see \cite{DaPZ,PrvtRcknr,Chw})  that the unique mild
solution is given by 
\begin{eqnarray}
 \label{eq1}
  X(t)=S(t)X_{0}+\int_{0}^{t}S(t-s)F(X(s))ds +O(t)
\end{eqnarray}
with the stochastic process $O$ given by the stochastic convolution 
\begin{eqnarray}
\label{OU}
O(t)=\int_{0}^{t}S(t-s)dW(s).
\end{eqnarray}
%
The study of numerical solutions of SPDEs is an active area of 
research and there is a growing literature on numerical methods for 
SPDEs (see \cite{allen98:_finit,Lakkis,KssrsZrrs,GTambueexpoM,Jentzen1,Jentzen2,Jentzen3,Jentzen4} and reference therein).

Our  numerical scheme is built on recent work by Jentzen and co-workers
\cite{Jentzen1,Jentzen2,Jentzen3,Jentzen4} that uses Taylor 
expansion and linear functionals of the noise for a spectral
Fourier--Galerkin discretisations of \eqref{adr} and obtained high order schemes in time. 
Let us  describe briefly these schemes.
Let $P_{N}$, $N\in \mathbb{N}$ be the spectral projection defined for
$u \in L^{2}(\Omega)$ by 
\begin{eqnarray} 
\label{sp}
P_{N}u=\sum_{i\in \mathcal{I}_{N}}(e_{i},u) e_{i},\qquad \qquad \mathcal{I}_{N}=\left\lbrace 1,2,...,N \right\rbrace^{d}.
\end{eqnarray}
The spectral Galerkin discretisation of \eqref{adr} yields the
following semi-discrete form 
\begin{eqnarray}
  \label{discretespect}
  dX^{N}=(A_{N}X^{N} +F_{N}(X^{N}))dt + d W^{N},
\end{eqnarray}
with $A_{N}=P_{N}A,\; F_{N}=P_{N} F$ and $W^{N}=P_{N}W$ and is a
diagonal system to solve for each Fourier mode.
For time stepping we make use of the standard $\varphi-$functions
\begin{eqnarray}
\varphi_{0}(\Delta t A_N)&=&e^{\Delta t A_{N}}\\ 
\varphi_{1}(\Delta t A_N)&=&(\Delta t\, A_{N})^{-1}\left( e^{\Delta
    t A_{N}}-I\right)= \frac{1}{\Delta t}\int_{0}^{\Delta t}
e^{(\Delta t- s)A_{N}}ds.
\end{eqnarray}

Jentzen and co-workers \cite{Jentzen3,Jentzen4} examine the following
two high order time stepping schemes which overcome the order barrier
(see \cite{Jentzen3}) of numerical schemes approximating \eqref{adr} 
\begin{equation}
  X_{m+1}^{N}=e^{\Delta t A_{N}}X_{m}^{N}+\Delta t \varphi_{1}(\Delta
  t A_{N}) F_{N}(X_{m}^{N})+P_{N} O_m
  \label{jentzen1}
\end{equation}
and 
\begin{equation}
 \label{jentzen2}
 Y_{m+1}^{N}=\varphi_{0}(\Delta t A_N) \left(Y_{m}^{N}+ \Delta t
 F_{N}(Y_{m}^{N}) \right)+ P_{N} O_m. 
\end{equation}
The process 
\begin{equation}
  O_{m} =\int_{t_{m}}^{t_{m+1}} e^{\left(t_{m+1}-s\right)A} dW
\end{equation}
has the exact variance in each Fourier mode as an Ornstein--Uhlenbeck process.
More  precisely, by assuming that the linear operator $A$ and the
covariance operator $Q$ have the same eigenbasis, applying the Ito
isometry in each mode yields 
\begin{eqnarray}
 ( e_{i},O_{m})=  e^{-\lambda_{i} \Delta
   t}\left(\dfrac{q_{i}}{2 \lambda_{i}} \left(1-e^{-2
       \lambda_{i}\Delta t} \right)\right)^{1/2}R_{i,m},
\label{eq:noiseupdate}
\end{eqnarray}
$i\in \mathcal{I}_{N}=\left\lbrace 1,2,3,...,N\right\rbrace^{2}$,
 $m=0,1,2..., M-1$ and $R_{i,m}$ are independent, standard normally
distributed random variables with means $0$ and variance $1$.  
In equation \eqref{eq:noiseupdate} the noise is termed to be computed using
a linear functional. 

Although schemes \eqref{jentzen1}-\eqref{jentzen2} are of higher order
in time, these improved convergence rates were only established under
seriously restrictive commutativity assumptions which exclude most
nonlinear Nemytskii operators.  This was recently overcome in \cite{xia}.
Another drawback is that to implement the schemes, the eigenfunctions
of the linear  operator $A$ and of the covarance operator $Q$
must coincide and furthermore must be known explicitly (see
\eqref{eq:noiseupdate}).  
To illustrate that this can be overcome with our spatial discretisation 
we solve the SPDE 
\begin{equation}
\label{couple}
  dX=\left(D \varDelta X -  \nabla \cdot (\textbf{q}
    X) + R(X)\right)dt+ dW,
\end{equation}
on a rectangular domain with mixed boundary conditions without
requiring information on the 
eigenvalues and eigenfunctions of the corresponding linear operator.
The velocity $\mathbf{q}$ in \eqref{couple} is obtained from the  following  steady state mass conservation equation and Darcy's law 
\begin{equation}
  \label{couple1}
  \nabla \cdot\mathbf{q} =0, \qquad \mathbf{q}=-\dfrac{\mathbf{k}}{\mu} \nabla p,
\end{equation}
where $\mathbf{k}$ is the heterogeneous permeability tensor, $p$ is
the pressure and $\mu$ is the dynamic viscosity of the fluid
\cite{sebastianb}. In \eqref{couple}, $R$ is the reaction function
which may be a Langmuir adsorption term which is globally Lipschitz
\cite{sebastianb} and $D>0$ the diffusion coefficient. 
Typically \eqref{couple}-\eqref{couple1} is solved using finite elements or finite
volumes as spectral Galerkin approach  is  not infeasible   due to the heterogeneous nature of the
permeability and  the fact that such problems often naturally give rise to non-uniform.
Our work differs from other finite element discretisations
\cite{allen98:_finit,Lakkis,KssrsZrrs,GTambueexpoM} where the noise
is considered directly in the finite element space.
We follow more closely \cite{Yn:04,Yn:05,Haus1} and introduce a
projection onto a finite number modes and a projection onto the finite
element space. The aim is to gain the flexibility of the finite
element (finite volume) discretisation to deal  with flow
and transport problems \eqref{couple}-\eqref{couple1}, complex
boundaries, mixed boundary conditions and  inhomogeneous boundary
conditions as well as reaching high order in time as in \cite{Jentzen3,Jentzen4}.

The paper is organised as follows. In \secref{sec:NumericalScheme} we
present the numerical scheme and assumptions that we make on the linear operator,
nonlinearity and the noise. 
We consider fairly weak conditions on nonlinear function $F$ as recently 
considered in \cite{xia}.
We then state and discuss our main results. These are convergence in the 
root mean square $L^{2}(\Omega)$ norm for reaction-diffusion equations and  advection--reaction--diffusion for 
spatially regular noise.
We present simulations in \secref{sect4}, these are applied both to a linear
example where we can compute an exact solution as well as a more
realistic model coming from model of the advection and diffusion of a
solute in a porous media with a non-linear reaction term.
We also show that, equipped with the eigenvalues and eigenfunctions of the operator
$\varDelta$ with Neumann or Dirichlet boundary conditions, we can apply
the new scheme with mixed boundary conditions 
without  explicitly having the eigenvalues and eigenfunctions for this case.
We present numerical results both for finite element and finite volume discretisations in space.
Finally, in \secref{sec:th1} and \secref{sec:th2}, 
 we present the proofs of the convergence theorems for the finite element discretisation. 
\section{Setting and Assumptions}
Let us start by presenting briefly the notation for the main function
spaces and norms that we use in the paper. 
We denote by $\Vert \cdot \Vert$ the norm associated to
the inner product $(\cdot ,\cdot )$ of the $\mathbb{R}-$Hilbert space $H=L^{2}(\Omega)$. 
For a Banach space $\mathcal{V}$ we denote  by $\Vert  \cdot \Vert_{\mathcal{V}}$ the norm of the space $\mathcal{V}$,
$L(\mathcal{V})$  the set of bounded linear mapping  from
$\mathcal{V}$ to $\mathcal{V}$, $L^{(2)}(\mathcal{V})$ the set of
bounded  bilinear mapping from  $\mathcal{V} \times \mathcal{V}$ to
$\mathbb{R}$ and $L_{2}(\mathbb{D},\mathcal{V})$ the Hilbert space of all equivalence 
classes of square integrable $\mathcal{V}-$valued random variables.

Let $Q: H\rightarrow H$ be a positive self adjoint operator, we  consider throughout this work the $Q$-Wiener process.
We denote the space of Hilbert--Schmidt operators from
$Q^{1/2}(H)$ to $H$   by $L_{2}^{0}:=\mathcal{L}_{2}(Q^{1/2}(H),H)= HS(Q^{1/2}(H),H)$
and the corresponding norm $\Vert . \Vert_{L_{2}^{0}}$ by
\begin{eqnarray*}
 \Vert l\,\Vert_{L_{2}^{0}} := \Vert l
 Q^{1/2}\Vert_{\mathcal{L}_{2}(H)}=\left( \underset{i \in \mathbb{N}^{d}}{\sum}\Vert
   l Q^{1/2} e_{i} \Vert^{2}\right)^{1/2},\qquad \ l\in L_{2}^{0}.
\end{eqnarray*}

Let  $\varphi : [0,T] \times \Omega \rightarrow L_{2}^{0} $  be a  $L_{2}^{0}-$valued 
predictable stochastic process with 
$\mathbb{P} \left[ \int_{0}^{t}\Vert \varphi \Vert_{L_{2}^{0}}^{2}ds < \infty \right]=1,\, t\in [0,T]$.
 We have the following
equality known as  the Ito's isometry
\begin{eqnarray*}
\mathbf{E} \Vert \int_{0}^{t}\varphi dW \Vert^{2}=\int_{0}^{t}
\mathbf{E} \Vert \varphi \Vert_{L_{2}^{0}}^{2}ds=\int_{0}^{t}
\mathbf{E} \Vert \varphi Q^{1/2} \Vert_{HS}^{2}ds. 
\end{eqnarray*}

Throughout the paper we assume that $\Omega$ is bounded and has a smooth boundary or is a convex polygon.
For convenience of presentation we take $A$ to be a self adjoint second order 
operator as this simplifies the convergence proof.
More precisely  
\begin{eqnarray}
\label{A}
A = \nabla \cdot \textbf{D} \nabla (.)+D_{0,0} \mathbf{I}
=\underset{i,j=1}{\sum^{d}}\dfrac{\partial }{\partial x_{i}}\left( D_{i,j}\dfrac{\partial 
 }{\partial x_{j}}\right)+D_{0,0} \mathbf{I},
\end{eqnarray}
where we assume that  $D_{i,j}=D_{j,i} \in L^{\infty}(\Omega)$ and that there exists a positive constant $c_{1}>0$ such that
 \begin{eqnarray}
\label{ellip}
\underset{i,j=1}{\sum^{d}}D_{i,j}(x)\xi_{i}\xi_{j}\geq c_{1}\vert \xi \vert^{2}  \;\;\;\;\;\;\forall \xi \in \mathbb{R}^{d}\;\;\; x \in \overline{\Omega}\;\;\; c_{1}>0.
\end{eqnarray}
The derivatives in \eqref{A} are understood in the sense of distributions (weak sense).
We introduce two spaces $\HH$ and $V$ where $\HH\subset V$. These spaces depend
on the choice of the boundary conditions and on the variational form associated to the operator
$A$. 
For Dirichlet boundary conditions we let 
\begin{eqnarray*}
V= \HH= H_{0}^{1}(\Omega)=\{v\in H^{1}(\Omega): v=0\quad
\text{on}\quad \partial \Omega\}.
\end{eqnarray*}
For Robin boundary conditions (Neumann boundary condition being a particular case)
we let $V=  H^{1}(\Omega)$ and
\begin{eqnarray*}
\HH = \left\lbrace v\in H^{2}(\Omega): \partial v/\partial
  \nu_{A}+\sigma v=0\quad \text{on}\quad \partial
  \Omega\right\rbrace, \qquad \sigma \in \mathbb{R}. 
 \end{eqnarray*}
 Note that $\partial v/\partial  \nu_{\mathcal{A}}$ is the normal derivative of $v$ and $ \nu_{\mathcal{A}}$ is the exterior pointing normal $\mathbf{n}=(n_i)$ 
 to the boundary of $\Lambda$ given  by 
\begin{eqnarray}
  \partial v/\partial
  \nu_{\mathcal{A}}= \underset{i,j=1}{\sum^{d}} n_i(x)D_{i,j}(x)\dfrac{\partial v}{ \partial x_j}.
\end{eqnarray}
Let $-A : \mathcal{D}(-A) \subset H \rightarrow H$ be the unbounded operator with domain  $\mathcal{D}(-A)$.
Under condition \eqref{ellip}, it is well known (see \cite{lions})
that the linear operator $-A$ generates an analytic semigroup $S(t):=e^{tA}$. 
Functions in $\HH$ can satisfy the boundary conditions. With the space $\HH$
in hand we can characterize the domain of the operator $(-A)^{r/2}$
and have the following norm equivalence \cite{lions,Stig,ElliottLarsson} 
for $r\in \{1,2\}$
\begin{eqnarray*}
C_{1}\Vert v \Vert_{H^{r}(\Omega)}&\leq & \Vert (-A)^{r/2} v\Vert \leq  C_{2}\Vert v \Vert_{H^{r}(\Omega)} 
\,\,\qquad \forall \,\, v\in \mathcal{D}((-A)^{r/2}),\\
\mathcal{D}((-A)^{r/2})&=&  \HH\cap H^{r}(\Omega)\,\,\qquad \qquad \qquad\quad \text{ (Dirichlet boundary conditions)},\\
\mathcal{D}((-A))&=&\HH,\,\,\,\mathcal{D}((-A)^{1/2})= H^{1}(\Omega),\qquad \,\,\text{(Robin boundary conditions)},
\end{eqnarray*} 
where $C_{1}$ and $C_{2}$ are positive constants. In fact  for Dirichlet, Robin and mixed boundary conditions
we have $V=\mathcal{D}((-A)^{1/2})$.
In the Banach space  $\mathcal{D}((-A)^{\alpha/2})$, $\alpha \in \mathbb{R}$, 
we will use the notation $ \Vert .\Vert_{\alpha} := \Vert (-A)^{\alpha/2}. \Vert $.

For our rigorous  convergence proof we  make the following assumptions on the linear operator $A$.
\begin{assumption}
  \label{assumption1} 
\textbf{[Linear operator]}
  The linear operator $-A$  given in  \eqref{A} is positive
  definite so there exists sequences of  positive 
  real eigenvalues $\{\lambda_{n}\}_{n\in \mathbb{N}^{d}}$ with $ \underset{i \in \mathbb{N}^{d}}{\inf}\lambda_{i} > 0$
  and an orthonormal basis in $H$ of eigenfunctions 
$\{e_{i}\}_{i \in\mathbb{N}^{d} }$ such that
  \begin{eqnarray*}
    -Av=\underset{i \in \mathbb{N}^{d}}{\sum}\lambda_{i}(
    e_{i},v)e_{i},\qquad \forall \quad v \in  \mathcal{D}(-A)  
  \end{eqnarray*}
  where $\mathcal{D}(-A)=\{ v \in H :\underset{ i\in \mathbb{N}^{d}}{\sum}\lambda_{i}^{2} |( e_{i},v) |^{2}<
  \infty\}$.
\end{assumption}
However we will show in a concrete example that for constant diffusion coefficient ($\mathbf{D}= D\, \mathbf{I},\, D>0$) with mixed boundary condition on rectangular grid, 
our scheme will be implemented with the well known eigenfunctions of  Laplace operator $\varDelta$ with Dirichlet or Neumann Boundary conditions. 
This flexibility can only be done  if  non-diagonal methods (finite element methods, finite volume method  and finite difference method) are used for space discretisation.

We recall some basic properties of the semi group $S(t)$ generated by $A$.
\begin{proposition}
\textbf{[Smoothing properties of the semi group\cite{Henry}]}\\
\label{prop1}
 Let $ \alpha >0,\;\beta \geq 0 $ and $0 \leq \gamma \leq 1$, then  there exist  $C>0$ such that
\begin{eqnarray*}
\Vert (-A)^{\beta}S(t)\Vert_{L(L^{2}(\Omega))} &\leq& C t^{-\beta}\qquad \text {for }\quad t>0\\
\Vert (-A)^{-\gamma}( \text{I}-S(t))\Vert_{L(L^{2}(\Omega))} &\leq& C t^{\gamma} \qquad \text {for }\quad t\geq 0.
\end{eqnarray*}
In addition,
\begin{eqnarray*}
(-A)^{\beta}S(t)&=& S(t)(-A)^{\beta}\quad \text{on}\quad \mathcal{D}((-A)^{\beta} )\\
\text{If}\quad \beta &\geq& \gamma \quad \text{then}\quad
\mathcal{D}((-A)^{\beta} )\subset \mathcal{D}((-A)^{\gamma} ),\\
\Vert D_{t}^{l}S(t)v\Vert_{\beta}&\leq& C t^{-l-(\beta-\alpha)/2} \,\Vert v\Vert_{\alpha},\;\; t>0,\;v\in  \mathcal{D}((-A)^{\alpha/2})\;\; l=0,1, 
\end{eqnarray*}
where $ D_{t}^{l}:=\dfrac{d^{l}}{d t^{l}}$.
\end{proposition}
The following assumption was recently used in \cite{xia}
and allows for more general $F$ than originally considered in \cite{Jentzen4}.
\begin{assumption}\label{assumption4}
\textbf{[Assumption on  nonlinear function $F$, and $Q$]}
For the noise, we assume that the covariance operator $Q$ satisfies
\begin{equation}
\label{noisenew}
\|(-A)^{\frac{\beta-1}{2}}Q^{\frac{1}{2}}\|_{\mathcal{L}_2(H)}<\infty,\text{ for some }\beta\in (0,2].
\end{equation}

For nonlinear function $F$, we assume that there exists a positive constant $L>0$ such that
$F$ satisfies either (a) or (b) below.

(a) $F:H\rightarrow H$ is Lipschitz, twice continuously
differentiable and satisfies for $X$, $Y$, $Z\in H$ 
\begin{eqnarray*}\label{asseq1}
 \|F(X)\|&\leq& L(1+\|X\|),  \qquad \|F'(Z)(X)\| \leq  L\|X\| \\
 \|(-A)^{-\eta/2}F''(Z)(X,Y)\|&\leq& L\|X\|\|Y\| \quad
  \text{for some }\eta \in[1,2).\label{asseq2} 
 \end{eqnarray*}
 Further for $X\in H$, $Z\in \mathcal{D}((-A)^{
   \frac{\min(\beta,1)}{2}})$, $\delta \in [1,2)$ 
$$
 \|(-A)^{-\frac{\delta}{2}}F'(Z)(X)\|\leq L(1+\|Z\|_{\min(\beta,1)})\|X\|_{-\min(\beta,1)}.$$
 
(b) $ F$ satisfies the following globally Lipschitz condition
$(H^{1}(\Omega),\Vert. \Vert_{H^{-1}(\Omega)})$
\begin{eqnarray*}
  \Vert F(Z)- F(Y)\Vert_{-1} \leq L \Vert Z- Y\Vert \qquad
\forall Z, Y  \in H. 
\end{eqnarray*}
\end{assumption}

\begin{remark}
 Using \cite[Lemma 2.3]{xia} we can easily see that if
 \eqref{noisenew} is satisfied then
$$\mathbf{E}  \Vert O(t)\Vert_\beta^2=\ \int_0^t \Vert (-A)^{\beta/2} S(t-s)\Vert^2 ds \leq C  \|(-A)^{\frac{\beta-1}{2}}Q^{\frac{1}{2}}\|_{\mathcal{L}_2(H)} \leq \infty.$$
\end{remark}
\begin{theorem}
  \label{existth}
  \textbf{[Existence, uniqueness \cite{DaPZ,Chw,PrvtRcknr,LrdPwllShrdlw}
and regularity \cite{GTambueexpoM,xia,kruse}.}
  Assume that the initial solution $X_{0}$ is an $F_{0}-$measurable
  $H-$valued random variable, the linear operator $A : \mathcal{D}(A)
  \subset H \rightarrow H$ is the generator of an analytic semigroup
  $S(t)=e^{t A}, \quad t\geq 0$, the noise $W$ is trace class and the
  nonlinear function $F: H\rightarrow  D \supseteq H$ is globally
  Lipschitz, where $D$ is a Banach space.
  There exists a mild solution $X$ to \eqref{adr} unique, up to
  equivalence among the processes, satisfying 
  \begin{eqnarray}
    \label{eq1.0}
    X(t) &=& S(t)X_{0}+\int_{0}^{t}S(t-s)F(X(s))ds+O(t)\\
    &=&\overline{X}(t)+O(t),
  \end{eqnarray}
  where  $O$ is the stochastic process  given by the stochastic
  convolution in \eqref{OU}.

  Furthermore if  \assref{assumption4} is satisfied with  the  corresponding $\beta$ in  \eqref{noisenew}, and $X_0\in L^p(\mathbb{D}, \mathcal{D}((-A)^{\beta/2})),\, p\in [2,\infty)$,
  then there exists a positive constant $C$ such that
  \begin{eqnarray}
   \underset{0\leq s\leq T}{\sup} \Vert X(t)\Vert_{L^p(\mathbb{D}, \mathcal{D}((-A)^{\beta/2}))} \leq C(1+\Vert X_0\Vert_{L^p(\mathbb{D}, \mathcal{D}((-A)^{\beta/2}))}).
  \end{eqnarray}
\end{theorem}
\section{Numerical scheme and main results}
\label{sec:NumericalScheme}
\subsection{Numerical scheme}
We consider discretisation of the spatial domain by a finite element 
triangulation.
Let $\mathcal{T}_{h}$ be a set of disjoint intervals  of $\Omega$ 
(for $d=1$), a triangulation of $\Omega$ (for $d=2$) or a set of
tetrahedra (for $d=3$) satisfying the standard regularity assumptions (see \cite{lions}).
Let $V_{h}\subset V $ denote the space of continuous functions that are
piecewise linear over the triangulation $\mathcal{T}_{h}$. 
To discretise in space we introduce two projections. Our first
projection operator 
$P_h$ is the $L^{2}(\Omega)$ projection onto $V_{h}$ defined  for  $u\in L^{2}(\Omega)$ by 
\begin{eqnarray}
(P_{h}u,\chi)=(u,\chi)\qquad \forall\;\chi \in V_{h}.
\end{eqnarray}
Then $A_{h}: V_{h}\rightarrow V_{h}$ is the discrete analogue of $A$ defined by
\begin{eqnarray}
( A_{h}\varphi,\chi)=a(\varphi,\chi)\qquad \varphi,\chi \in V_{h},
\end{eqnarray}
where $a(\ ,\ )$ is the corresponding bilinear form associated to the operator $A$.
We denote by $S_h$ the semigroup generated by the operator $A_h$.
The second projection $P_{N}$, $N \in \mathbb{N} $ is the projection
onto a finite number of spectral modes ${e_i}$ defined for $u\in L^{2}(\Omega)$ by  
$$P_{N}u=\sum_{i\in \mathcal{I}_{N}} (e_{i},u) e_{i},$$
where  $\mathcal{I}_{N}=\left\lbrace 1,2,...,N \right\rbrace^{d}$.

The semi--discrete version of the problem (\ref{adr}) is to
find the process $X_{h}(t)=X_{h}(.,t) \in V_{h}$ such  that for $t\in[0, T]$,
\begin{eqnarray}
\label{dadr}
dX_{h}=(A_{h}X_{h} +P_{h}F(X_{h}))dt + P_{h} P_N d W,\qquad
 X_{h}(0)=P_{h}X_{0}.
\end{eqnarray}
The solution of \eqref{dadr} is given by 
$$
X_h(t)=S_h(t)X_h(0)+\int_{0}^{t}S_h(t-s)F(X_h(s))ds+\int_{0}^{t}S_h(t-s)P_{h} P_N dW.
$$
Set $O_{h}(t)$ and $O^{h}(t)$ two $V_{h}$-valued  stochastic convolutions  defined by
\begin{eqnarray}
\label{dconv}
 O_{h}(t)=\int_{0}^{t}S_h(t-s)P_{h} P_N dW,\qquad O^{h}(t)=P_{h} P_N \int_{0}^{t}S(t-s)dW.
\end{eqnarray}
In order to build our scheme based on semi-implicit discretisation in time and linear 
 functional of the noise, we used  the approximation $O_{h}(t)\approx O^{h}(t)$.
Notice that the two stochastic convolutions are the space approximations of  the convolution $O(t)$ defined in \eqref{OU}. 
 We therefore have the following  semi-discrete solution
$$
X^h(t)=S_h(t)X^h(0)+\int_{0}^{t}S_h(t-s)F(X^h(s))ds+P_{h} P_N \int_{0}^{t}S(t-s)dW.
$$
 We denote by $\overline{X}^{h}$ the solution of the random system 
$$
\overline{X}^h(t)=S_h(t)X^h(0)+\int_{0}^{t}S_h(t-s)F(X^h(s))ds.
$$
As in \eqref{eq1}, by splitting we have
$$
X^{h}(t)=\overline{X}^h(t)+  P_{h} P_N O(t).
$$
We now discretise in time by a semi--implicit method to get
the  fully discrete approximation of $\overline{X}^{h}$ defined
by $Z_{m}^{h}$ 
\begin{eqnarray}
  Z_{m}^{h}= \Shdt^m P_{h} X_{0}+\Delta\,t\,\underset{k=0}{\sum^{m-1}}
\Shdt^{(m-k)}P_{h}F( Z_{k}^{h}+P_{h}P_{N}O(t_{k})),
\end{eqnarray}
where
\begin{equation}
  S_{h,\Delta t}:=(\text{I}-\Delta t\;A_{h})^{-1}.
\label{eq:Shdt}
\end{equation}
It is straightforward to show that 
\begin{eqnarray}
  Z_{m+1}^{h}=\Shdt 
\left(Z_{m}^{h} +\Delta t P_{h}F(Z_{m}^{h}+P_{h}P_{N}O(t_{m}))\right).  
\label{eq:Zmh}
\end{eqnarray}
Finally we can define our approximation $X_{m}^{h}$ to $X(t_{m})$, the solution of equation (\ref{adr}) by
\begin{eqnarray}
 X_{m}^{h}=Z_{m}^{h}+P_{h}P_{N}O(t_{m}).
\label{eq:XbarO}
\end{eqnarray}
Therefore
\begin{eqnarray}
  \label{new}
 X_{m+1}^{h}= \Shdt 
\left(X_{m}^{h} +\Delta t \,
P_{h}F(X_{m}^{h})-P_{h}P_{N}O(t_{m})\right) +P_{h}P_{N}O(t_{m+1}),
\end{eqnarray}
where  according  to \eqref{OU}, we generate $O(t_{m+1})$ from  $O(t_{m})$  by
\begin{equation}
\label{eq:O}
O(t_{m+1}) = e^{\Dt A} O(t_{m}) + \int_{t_{m}}^{t_{m+1}}
e^{(t_{m+1}-s)A}dW(s). 
\end{equation}
%
The new modified scheme \eqref{new},\eqref{eq:O} uses a finite
element discretisation and projects the linear functional of the noise on the space $V_h$
and hence we expect superior  approximation properties over a standard
semi-implicit Euler--Maruyama discretisation for a finite element
discretisation (given in \eqref{standard} below). 

From the Sobolev embedding theorems (see for example \cite[Theorem 3.10]{EP})
 we formulate the following remark which allows us to replace
$ O^{h}(t)=P_h P_N O(t)$ by the interpolation of the convolution $ P_N O(t)$ at the finite element nodes.  Therefore to simulate  the
convolution $O^{h}(t)$ we can just evaluate the convolution $ P_N O(t)$ at the  nodes of the finite element mesh.
\begin{remark}
\label{evaluaten}
 If the noise is regular enough in space, i.e. $O(t)\in H^{k+1}(\Omega)$ ($k$ being determined such $H^{k+1}(\Omega) \subset C(\overline{\Omega})$), 
then the  $L^{2}(\Omega)$ orthogonal  projection $P_h$ in the stochastic convolution $ O^{h}(t)$ 
can be replaced by  the interpolation operator $I_h :H^{k+1}(\Omega) \rightarrow V_h$  defined for $u \in H^{k+1}(\Omega)$ by
\begin{eqnarray}
I_h(u)=\underset{i=1}{ \sum^{N_h}}u(a_{i})\varphi_{i},
\end{eqnarray}
where $a_1 ,..., a_{N_h}$ are the finite element nodes, $N_h=\text{dim}(V_h)$ and $\varphi_{1},...,\varphi_{N_h}$  the corresponding nodal basis with $\varphi_i(a_j)=\delta_{i,j}$.
\end{remark}
The standard semi-implicit Euler- Maruyama scheme  for (\ref{adr}) is given by
\begin{eqnarray}
  \label{standard}
 Y_{m+1}^{h}&=& \Shdt 
\left(X_{m}^{h} +\Delta t \,
P_{h}F(Y_{m}^{h})+P_{h} \Delta W_{m}^{N} \right)\\
\Delta W_{m}^{N}&:=& W_{t_{m+1}}^{N}-W_{t_{m}}^{N}=\sqrt{\Delta t} \underset{i \in \mathcal{I}_{N}}{\sum}\sqrt{q_{i}} R_{i,m}e_{i} \nonumber\\
\end{eqnarray}
where $R_{i,m}$ are independent, standard normally distributed random
variables with mean $0$ and variance $1$.  This standard  scheme will be used in \secref{sect4} for comparison with the new scheme.
 We use the Monte Carlo method to approximate  the discrete root mean square $L^{2}$ norm of the error on a regular 
mesh with size $h$ at the final time $T=M\Dt$. Indeed we use that
\begin{eqnarray}
\label{error}
 \left(\mathbf{E}\Vert X(T)-\xi_{M}^{h}\Vert^{2}\right)^{1/2}&= &\left(\mathbf{E}\Vert X(.,T)-\xi_{M}^{h}(.)\Vert^{2}\right)^{1/2} \nonumber\\
&\approx&  \left(\dfrac{h^d}{K} \underset {l=1}{\sum ^{K}}  \underset {i=1}{\sum ^{N_h}} \left(X(a_{i},T)-\xi_{M}^{h}(a_{i})\right)^{2}\right)^{1/2},
\end{eqnarray}
 where $\xi_{M}^{h}$ is  either $X_{M}^{h}$ or $Y_{M}^{h}$ (the numerical solutions from the final step in \eqref{new} or \eqref{standard} 
for each sample $l$),  $K$ is the number of sample solutions and $X(T)$ 
is  the 'exact' solution for the sample $l$ that we will specify for each example in \secref{sect4}.
\subsection{Main results: strong convergence in $L^2$}
Throughout the article we let $N$ be the number of terms of truncated
noise, 
and let $t_m=m\Dt \in (0,T]$, where $T=M\Dt$ for $m,M\in\N$. We take $C$
to be a constant that may depend on $T$ and other parameters but not
on $\Dt$, $N$ or $h$.
We examine strong convergence for the two distinct assumptions on the
nonlinearity $F$ given in \assref{assumption4}. We present the
theorems here and the proofs may be found in \secref{sec:proofs}.
When the non-linearity satisfies the Lipschitz
condition of \assref{assumption4} (a) we have the following theorem.
\begin{theorem}
\label{th35}
Suppose that \assref{assumption1} holds, and  the noise covariance $Q$ and  the non-linearity $F$ satisfy respectively  \eqref{noisenew} and  condition  a)  of 
\assref{assumption4}.
Let  $X(t_m)$  be the mild  solution of equation \eqref{adr} represented
by (\ref{eq1})  and $ X_m^{h}$ be the numerical
approximations through scheme \eqref{new}.
Let $\beta \in(0,2)$ be as defined in \eqref{noisenew}. Then if $X_0
\in L(\mathbb{D},\mathcal{D}((-A)^{\beta/2}))$
\begin{eqnarray*}
 \lefteqn{\left(\mathbf{E}\Vert X(t_m)-X_{m}^{h}\Vert^{2}\right)^{1/2}}\\
 & \leq &C \left( h^{\beta}+\Delta t^{\min(\beta,1)} +  \Delta t \vert \ln (\Delta t)\vert+\left(\underset { j \in \mathbb{N}^{d} \backslash
       \mathcal{I}_{N}} {\inf}  \lambda_{j}\right)^{-\beta/2}\right).
\end{eqnarray*}
\end{theorem}
The convergence in the mean square  $L^2(\Omega)$ 
norm where the non-linearity satisfies the Lipschitz condition
from  $H^{-1}(\Omega)$ norm to $L^2(\Omega)$ (\assref{assumption4} (b)) is given in the following theorem.
\begin{theorem}
\label{th25}
Suppose that \assref{assumption1} holds, and  the noise covariance $Q$ and  the non-linearity $F$ satisfy respectively  \eqref{noisenew} and  condition  b)  of 
\assref{assumption4}.
Let  $X(t_m)$  be the mild  solution of equation \eqref{adr} represented
by (\ref{eq1})  and $ X_m^{h}$ be the numerical
approximations through scheme \eqref{new}.
Let $\beta \in(0,2)$ be as defined in \eqref{noisenew}, Then if $X_0
\in L(\mathbb{D},\mathcal{D}((-A)^{\beta/2}))$
\begin{eqnarray*}
  \lefteqn{\left(\mathbf{E}\Vert X(t_m)-X_{m}^{h}\Vert^{2}\right)^{1/2}}\\
  & \leq
 &C \left( h^{\beta}+\Delta t^{\min(\beta/2,1/2)} + \Delta t \vert \ln (\Delta t)\vert+ \left(\underset { j \in \mathbb{N}^{d} \backslash
       \mathcal{I}_{N}} {\inf}  \lambda_{j}\right)^{-\beta/2}\right).
\end{eqnarray*}
\end{theorem}
We note that Theorem \ref{th25} covers the case of
advection-diffusion-reaction SPDEs, such as that arising in our
example from porous media. However, we see a reduction in the
convergence rate compared to Theorem \ref{th35}.

If we denote by  $N_{h}$ the number of
vertices in the finite element mesh then it is well known (see for 
example \cite{Yn:05}) that  if $N \geq N_{h}\ $ then 
$$
  \left( \underset { j \in \mathbb{N}^{d} \backslash \mathcal{I}_{N}}
    {\inf} \lambda_{j}\right )^{-\beta/2} \leq C h^{\beta}.
$$
As a consequence the estimates in Theorem
\ref{th35}  and Theorem \ref{th25} can be expressed as functions of $h$
and $\Delta t $ only, and it is the error from the finite element
approximation that dominates. If $N\leq N_h$ then it is the error from
the projection $P_N$ of the noise onto a finite number of modes that
dominates.
\begin{remark}
Comparing our errors estimates with the one of standard semi-implicit Euler-Maruyama scheme \eqref{standard} given in \cite{kruse}, we observe that our new scheme  is higher
order in time (almost twice  the order of  the standard semi-implicit Euler-Maruyama scheme) when the nonlinear $F$ satisfies the \assref{assumption4} (a).  
\end{remark}

\section{Numerical Simulations}
\label{sect4}
We consider two example SPDEs for our numerical simulations.
Our first example is linear and we can construct an explicit solution.
We will examine to different types of noise for this case - one where
$A$ and $Q$ have the same eigenfunctions and one where they do not.
Our second example is motivated from realistic porous
media flow and has mixed boundary conditions and in this more
challenging example we assume the eigenfunctions of $A$ and $Q$ coincide.

In all cases the linear operator $A$ is linked to the Laplace
operator $\varDelta$  with homogeneous Neumann boundary conditions on
the domain 
$\Omega=[0,L_{1}]\times [0,L_{2}]$. 
The eigenfunctions $\{e_{i}^{(1)}\otimes e_{j}^{(2)}\}_{i,j\geq 0}
$ of the operator $-\varDelta$ here are given by 
\begin{eqnarray*}
e_{0}^{(l)}(x)=\sqrt{\dfrac{1}{L_{l}}},\qquad 
e_{i}^{(l)}(x)=\sqrt{\dfrac{2}{L_{l}}}\cos(\lambda_{i}^{(l)}x), \qquad \lambda_{0}^{(l)}=0,\qquad
\lambda_{i}^{(l)}=\dfrac{i \pi }{L_{l}}
\end{eqnarray*}
where $l \in \left\lbrace 1, 2 \right\rbrace,\, x\in \Omega$ and  $i=1, 2, 3, \cdots$
with the corresponding eigenvalues $ \{\lambda_{i,j}\}_{i,j\geq 0} $ given by 
$\lambda_{i,j}= (\lambda_{i}^{(1)})^{2}+ (\lambda_{j}^{(2)})^{2}$.

We use two types of noise in our simulations. In both examples,
we take the eigenvalues 
\begin{eqnarray}
\label{noise2}
 q_{i,j}=\left( i^{2}+j^{2}\right)^{-(\beta +\epsilon)}, \, \beta>0,
\end{eqnarray} 
in the representation \eqref{eq:W1} for some small $\epsilon>0$. Here
the noise and the linear operator have are supposed to have the same
eigenfunctions. 
We obviously have 
\begin{eqnarray*}
\underset{(i,j) \in \mathbb{N}^{2}}{\sum}\lambda_{i,j}^{\beta-1}q_{i,j}<  \pi^{2}\underset{(i,j) 
\in \mathbb{N}^{2}}{\sum} \left( i^{2}+j^{2}\right)^{-(1+\delta)} <\infty,
\end{eqnarray*}
thus \eqref{noisenew} in \assref{assumption4} is satisfied.  We take
$\beta\in \{1,2\}$ below.

To illustrate that our scheme can be used when $A$ and $Q$ have
different eigenfunctions we also consider for the first example an
exponential covariance for $Q$ so that 
 $$
\mathbf{E} W((x_1,y_1),t)W((x_2,y_2),t')=C_r((x_1,y_1),(x_2,y_2))\min(t,t')$$  with
\begin{eqnarray}
\label{expo}
 C_{r}((x_{1},y_{1});(x_{2},y_{2}))=\dfrac{\Gamma}{4 b_{1}b_{2}} \exp \left(-\dfrac{\pi}{4}\left[\dfrac{\left( x_{2}-x_{1}\right)^{2} }{c_{1}^{2}}+ \dfrac{\left( y_{2}-y_{1}\right)^{2} }{c_{2}^{2}}\right] \right)  
\end{eqnarray}
where $b_{1},b_{2}$ are  spatial correlation lengths in $x$ and $y$
and  $\Gamma>0$. Here the regularity of the solution depends of the
correlation lengths $b_{1},b_{2}$ (small values implying less
regularity). Similar to \cite{shardlow05,AtThesis} we can obtain an
approximation to the eigenvalues $q_{i,j}$ of $Q$.
\begin{proposition}
  \label{Prop:QEVALS}
  Assuming  that  $b_{i}\ll L_{i}$,  the projection of the noise on
  the eigenvectors of $A$ 
  yields  the following coefficients in the representation \eqref{eq:W1}
  \begin{eqnarray}
    \label{eq:Qexpevls}
    q_{i,j}= \Gamma \exp\left[ -\dfrac{1}{2 \pi}\left((\lambda_{i}^{(1)}b_{1})^{2}+(\lambda_{j}^{(2)}b_{2})^{2}\right) \right].
  \end{eqnarray}
\end{proposition}
The proof of this proposition can be found in
\secref{sec:Prop:QEVALS}. We see from \eqref{eq:Qexpevls} that again $Q$ is in trace class.
Below we take $b=b_1=b_2$ and use \eqref{eq:Qexpevls} in \eqref{eq:W1}. 


In all our simulations, the noise is truncated and we take $\vert \mathcal{I}_N\vert =N_{h}= \text{dim}(V_h),$ then $N\geq N_{h}^{1/d}$ 
 as suggested in \cite{Yn:05,newstig,newstigt} to avoid the reduction of the orders of convergence. 
In the case of the exponential covariance function  \eqref{expo}, work in \cite{newstigt} suggested that $\vert \mathcal{I}_N \vert $ can be $\ln(N_{h})$
 and the orders of convergence are still preserved. 
 
 In the implementation of our modified scheme
 \eqref{new},\eqref{eq:O} at every time step,
$O(t_{k+1})$ is generated using  $O(t_{k})$ from the following relation 
\begin{eqnarray*}
O(t_{k+1})= e^{\Delta t A} O(t_{k}) +\int_{t_{k} }^{t_{k+1}} e^{(t_{k+1}-s)A}dW(s),
\end{eqnarray*}
where $O(0)=0$.
We expand in Fourier space and apply the Ito isometry
in each mode and project onto $N$ modes to obtain for $k=1,2,\ldots,M-1$
\begin{eqnarray}
\label{singnoise}
( e_{i},O(t_{k+1}))= e^{-\lambda_{i}
 \Delta t}(e_{i},O(t_{k}))+ \left(\dfrac{q_{i}}{2 \lambda_{i}} \left(1-e^{-2
 \lambda_{i}\Delta t} \right)\right)^{1/2}R_{i,k},
\end{eqnarray}
where  $ R_{i,k}$ are independent, standard normally distributed random
variables with mean $0$ and variance $1$, and 
$i\in \mathcal{I}_{N}=\left\lbrace1,2,3,...,N\right\rbrace^{2}$. 
The noise is then projected onto the 
finite element space by $P_h$.  As we notice in \rmref{evaluaten}, for continuous noise $P_hP_N O(t) $ can be replaced by $I_hP_N O(t)$ i.e 
the evaluation of $P_N O(t)$ at the mesh  vertices. 
If the noise is not smooth then $P_h P_N O(t)$ is evaluated following the work in \cite[Section 5]{Yn:05}  for $P_h W$. Indeed, by setting 
$P_hP_N O(t)= \underset{i=1}{ \sum^{N_h}} \alpha_i^{1/2} \varphi_i $,  as $(e_{i},O(t_{k}))$  is known from \eqref{singnoise}, 
the coefficients $\alpha_i$  are found 
by solving the linear system
$$ \underset{i=1}{ \sum^{N_h}}(e_{i},O(t_{k}))^2(e_{i},\varphi_j)^{2} = \underset{i=1}{ \sum^{N_h}} \alpha_i( \varphi_i,\varphi_j)^{2},\,\,\,\,\,\,\,\qquad \qquad j=1,2,....,N_h, $$
 where $(\varphi_{i})_{1\leq i \leq N_{h}}$ is the nodal basis with $\varphi_i(a_j)=\delta_{i,j}$.  
 \subsection{A linear reaction--diffusion equation with exact solution}
As our first simple example we consider the reaction diffusion equation 
\begin{eqnarray}
\label{linear}
 dX=(D \varDelta X -0.5 X)dt+ dW, 
\qquad \text{given } \quad X(0)=X_{0},
\end{eqnarray}
on the time interval $[0,T]$.
We take $L_1=L_2=1$. 
Notice that $A=D \varDelta $ does not satisfy \assref{assumption1} as
$0$ is an eigenvalue. 
For simulation, one can eliminate the eigenvalue $0$, use the
perturbed operator $A=D \varDelta+\epsilon \mathbf{I},\, \epsilon>0$  
or eliminate the node with eigenvalue $0$ if $q_{0}=0$.  
The exact solution of \eqref{linear} is known. Indeed,
the decomposition of \eqref{linear} in each eigenvector node yields the following  Ornstein-Uhlenbeck process
\begin{eqnarray}
\label{exact0}
 dX_{i}=-(D \lambda_{i}+0.5)X_{i}dt+ \sqrt{q_{i}}d\beta_{i}(t)\qquad i \in \mathbb{N}^{2}.
\end{eqnarray}
This is a Gaussian process with the mild solution
\begin{eqnarray}
X_{i}(t)= e^{-k_{i}t}X_{i}(0)+  \sqrt{q_{i}}\int_{0}^{t}e^{k_{i}(s-t)} d \beta_{i}(s),\quad k_{i}=D \lambda_{i}+0.5.
\end{eqnarray}
Applying the Ito isometry yields the following exact variance of $X_{i}(t)$
\begin{eqnarray}
 \text{Var}(X_{i}(t))=\dfrac{q_{i}}{2 k_{i}}\left(1-e^{-2 \,k_{i}\,t}\right).
\end{eqnarray}
During simulation, we compute the exact solution recurrently as 
\begin{eqnarray}
\label{exact}
X_{i}^{m+1}&=& e^{-k_{i} \Delta t}X_{i}^m+  \sqrt{q_{i}}\int_{t_m}^{t_{m+1}}e^{k_{i}(s-t)} d \beta_{i}(s)\nonumber\\
            &=& e^{-k_{i} \Delta t }X_{i}^m +\left(\dfrac{q_{i}}{2 k_{i}}\left(1-e^{-2 \,k_{i}\,\Delta t}\right)\right)^{1/2}R_{i,m},
\end{eqnarray}
where $R_{i,m}$ are independent, standard normally distributed random
variables with mean $0$ and variance $1$.
The expression in \eqref{exact}  allows us to use the same set of
random numbers for both the exact and the numerical solutions.  

Our function $F(u)=-0.5 u $ is linear and obviously satisfies
\assref{assumption4} (a). In our simulation we take $D=1$ and
$X_{0}=0$. With $X_0=0$ the mild solution $X$ satisfies the regularity
required in Theorem \ref{th35}.
%
%

\begin{figure}[!ht]
  \subfigure[]{
    \label{FIGIa}
    \includegraphics[width=0.5\textwidth]{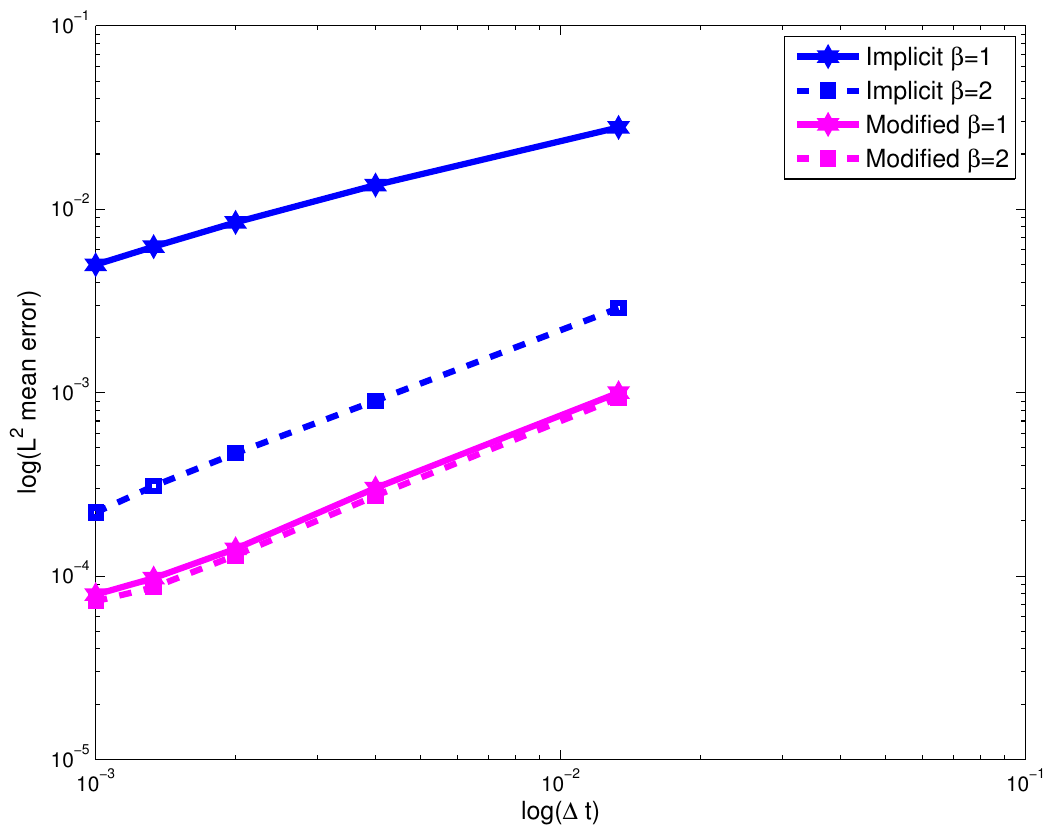}}
  \subfigure[]{
    \label{FIGIb}
    \includegraphics[width=0.5\textwidth]{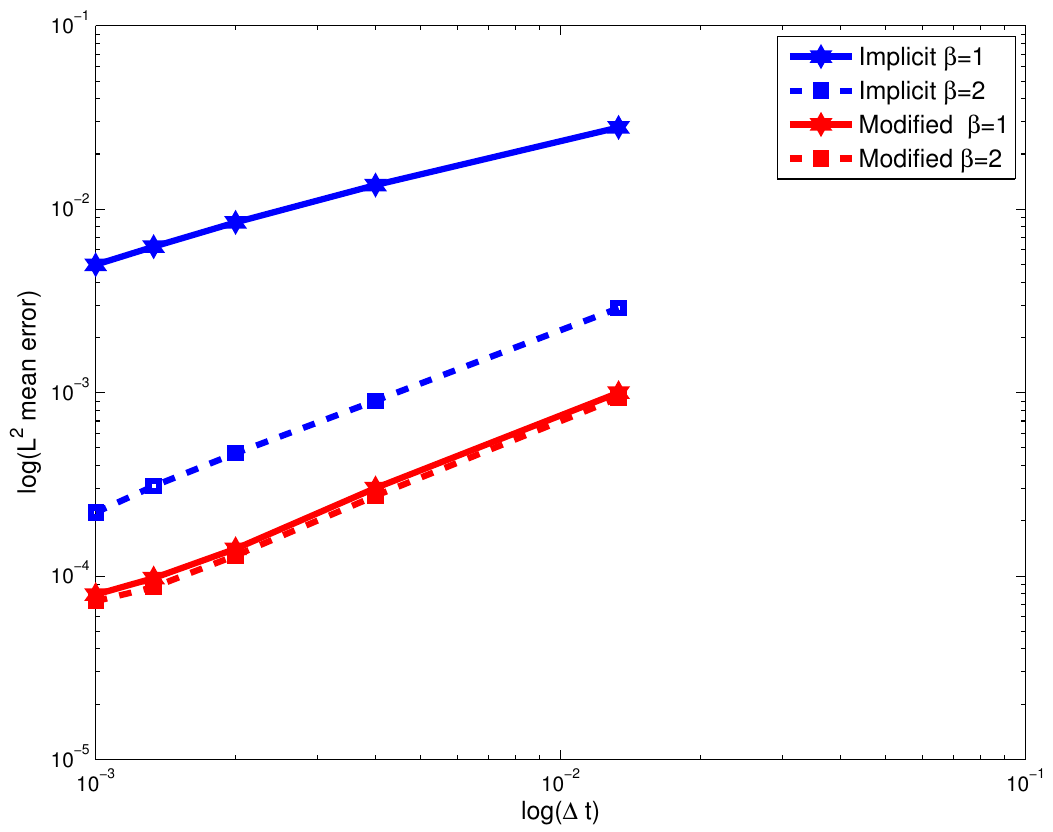}}
  \caption{Convergence in the root mean square $L^{2}$ norm at $T=1$ as a
    function of $\Dt$ for \eqref{linear}. We show convergence for noise  where 
    $\beta \in \{1,2\}$ and $\delta=0.05$ in relation \eqref{noise2}  for 
    finite element (a) and finite volume (b) space discretizations. We
    also show convergence of the standard semi--implicit scheme for the
    finite volume discretization. We used here 30 realizations. 
    Note that graphs with finite element and finite volume methods
    give the same errors. Note that $\beta=2$ represents here the case
    $\beta=2-\epsilon$, with positive $\epsilon$ small enough.} 
  \label{FIGI} 
\end{figure}

We examine both the finite element and the finite
volume discretization in space. For the cell center finite volume
discretization we take $h=\Delta x=\Delta y =1/100$. 
The finite element triangulation was constructed so that the center of the
control volume for the finite volume method was a vertex in finite element
mesh. To examine the error we use the exact solution \eqref{exact} 
and the discrete error estimate given in
\eqref{error}.
In \figref{FIGIa}, we see that the observed rate of convergence for
the finite element discretization agrees with \thmref{th35}.
The rate of convergence in $\Dt$ is very close to $1$ for $\beta\in
\{1,2\}$ in our new modified scheme.
We also see that due to the regularity of the mesh, the finite element
(\figref{FIGIa}) and finite volume methods (\figref{FIGIb}) give essentially the same errors.
More importantly we see that the new modified scheme is more accurate than the standard
implicit Euler--Maruyama scheme.
Indeed, we observe numerically a slower rate of convergence
of $0.65$ (for $\beta=1$) and $0.98$ (for $\beta=2$) of the standard scheme
compared to $ 0.9960$ (for $\beta=1$) and  $1.0074$ (for $\beta=2$) with the modified scheme.
We also observe that the error decreases as the regularity increases
from $\beta=1$ to $\beta=2$.  

\begin{figure}[!th]
  \begin{center}
    \includegraphics[width=1\textwidth,height=0.4\textheight]{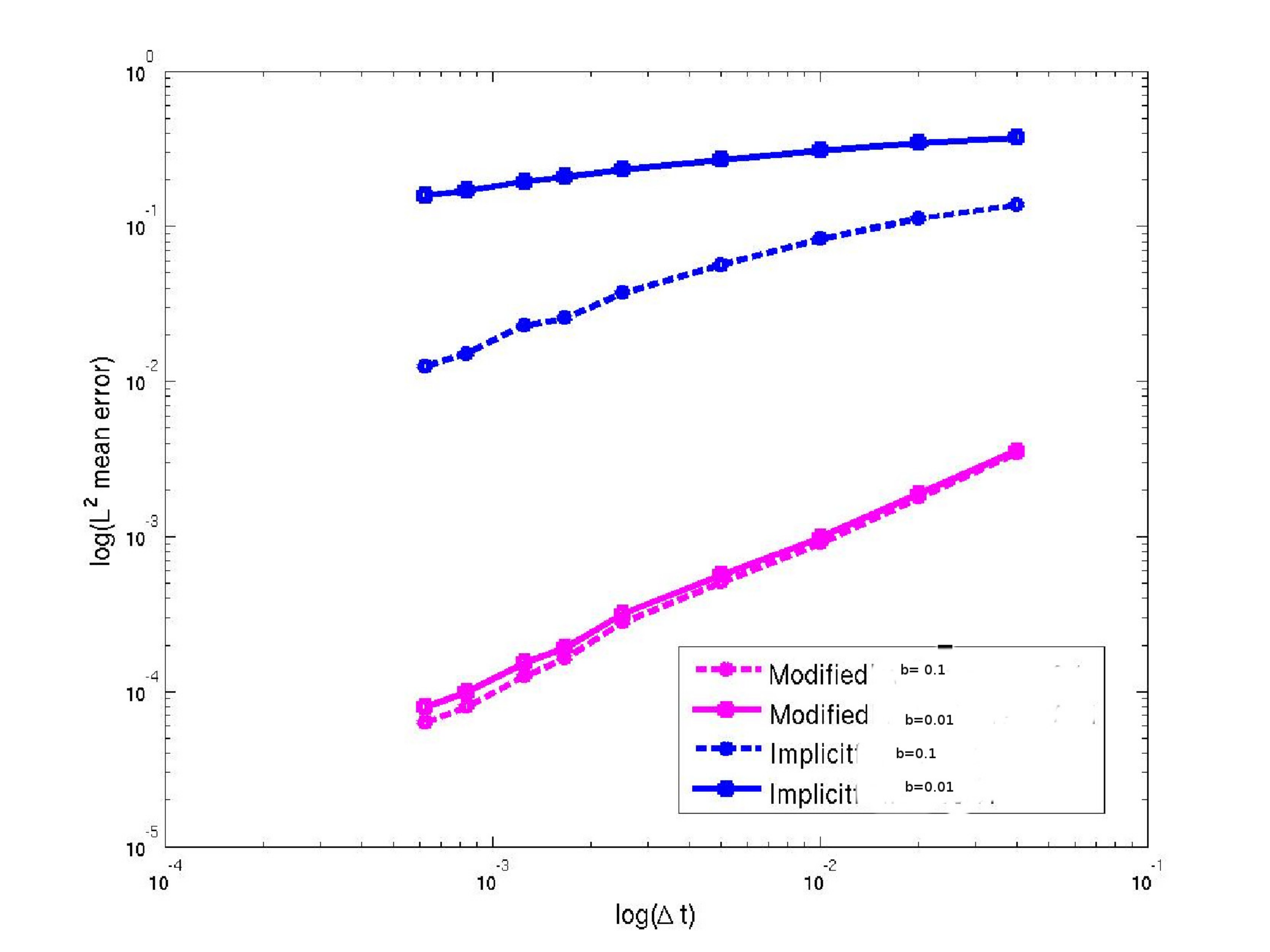}
    \caption{Convergence in the root mean square $L^{2}$ norm at $T=1$ as a
      function of $\Dt$ with the noise having exponential correlation
      function \eqref{expo} for a finite element discretisaton of
      \eqref{linear}. We used here 30 realizations and two 
      correlation length, $b=b_{1}=b_{2}=0.1$ and $b=b_{1}=b_{2}=0.01$,
      $\Gamma=1$.} 
    \label{FIGII} 
  \end{center}
\end{figure}

\figref{FIGII} shows convergence with exponential correlation
\eqref{expo}. Again the new
modified scheme is more accurate than the standard
semi-implicit Euler--Maruyama scheme. We observe numerically a slower rate of convergence
of $0.2098$ (for $b_{1}= b_{2}=0.01$) and $0.5918$ (for $b_{1}= b_{2}=0.1$) of the standard scheme
compared to $0.9145$ (for  $b_{1}= b_{2}=0.01$) and  $0.9618$ (for
$b_{1}= b_{2}=0.1$) with the modified scheme.  

\subsection{Stochastic advection diffusion reaction with mixed boundary conditions}
For our second, and more challenging example, we consider the stochastic advection
diffusion reaction SPDE \eqref{couple}, with $D=10^{-2}$
and mixed Neumann-Dirichlet boundary conditions on $\Omega=[0,1]\times[0,1]$. 
The Dirichlet boundary condition is $X=1$ at $\Gamma=\{ (x,y) :\; x =0\}$ and 
we use the homogeneous Neumann boundary conditions elsewhere.

Our goal here is to show that with the well known eigenvalues and eigenfunctions of the operator
$\varDelta$ with Neumann (or Dirichlet) boundary conditions, we can apply 
the new scheme to mixed boundary conditions for the operator $A=D\varDelta$ without explicitly having
eigenvalues and eigenfunctions of $A$. We also show that the modified scheme is more accurate than the 
standard semi-implicit Euler--Maruyama method. 

\begin{figure}[h!]
\subfigure[]{
\label{FIG022a}
\includegraphics[width=0.50\textwidth]{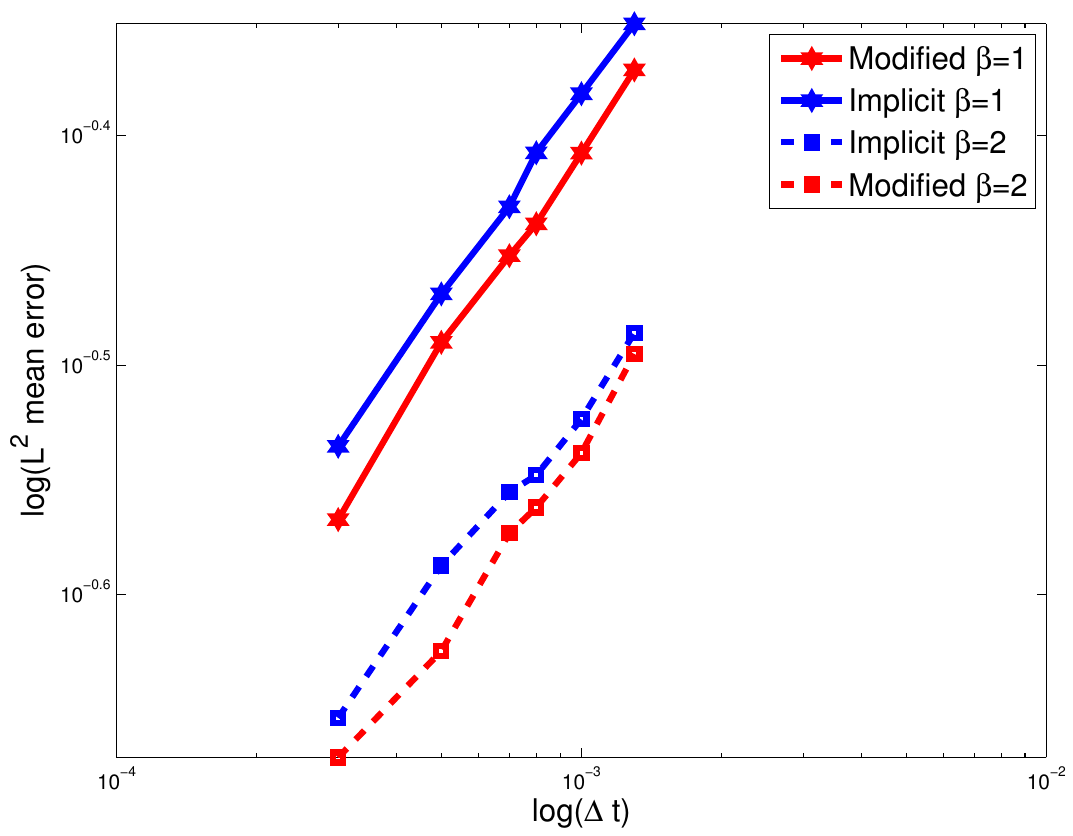}}
\subfigure[]{
\label{FIG022b}
 \includegraphics[width=0.50\textwidth]{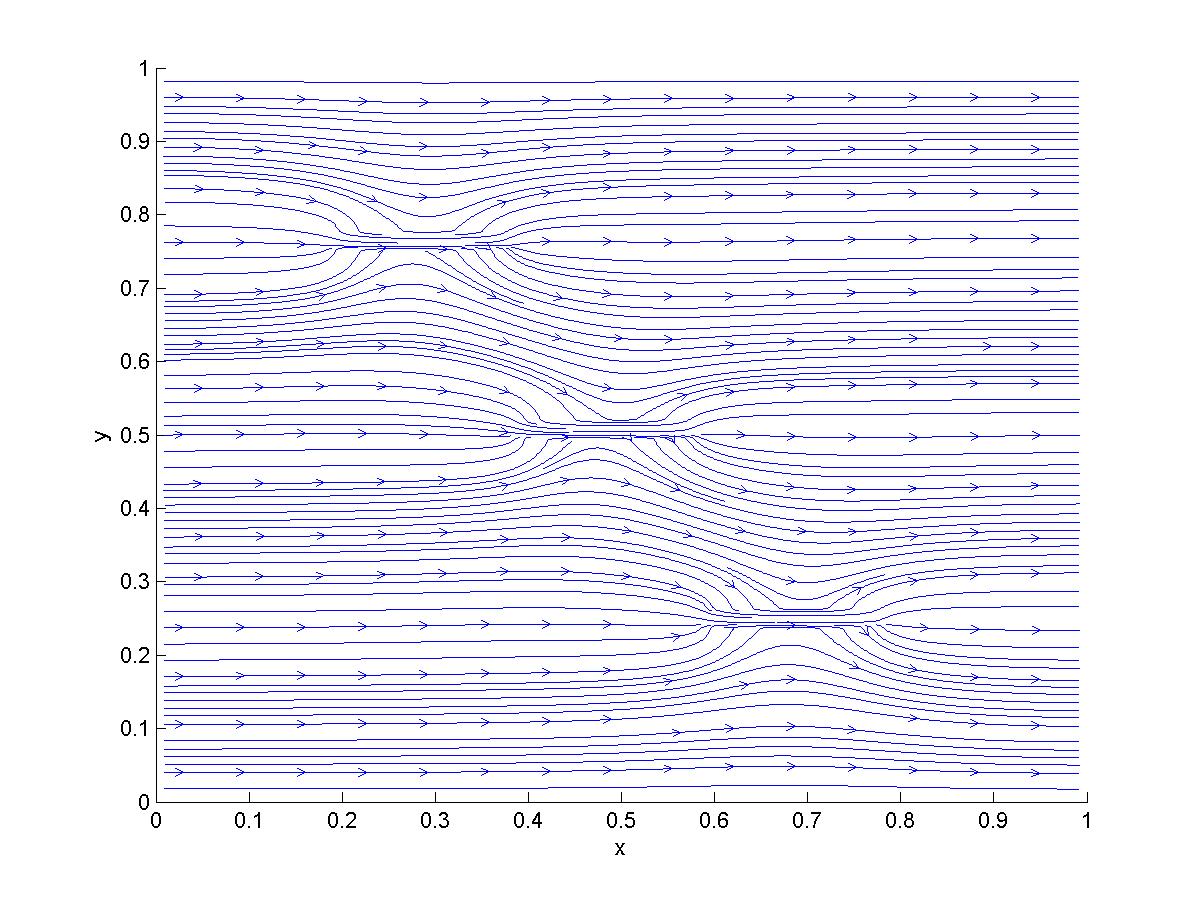}}
\subfigure[]{
 \label{FIG022c}
\includegraphics[width=0.50\textwidth]{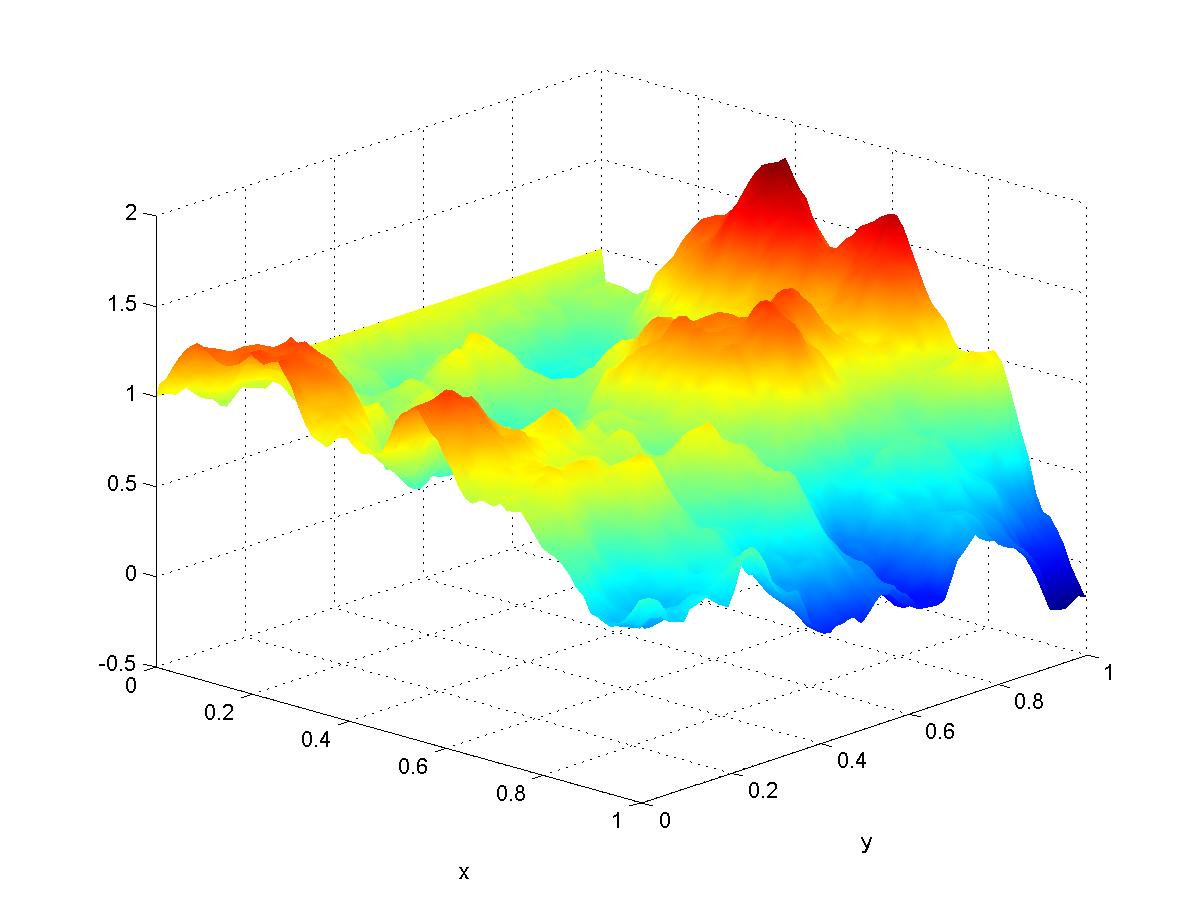}}
\subfigure[]{
\label{FIG022d}
\includegraphics[width=0.50\textwidth]{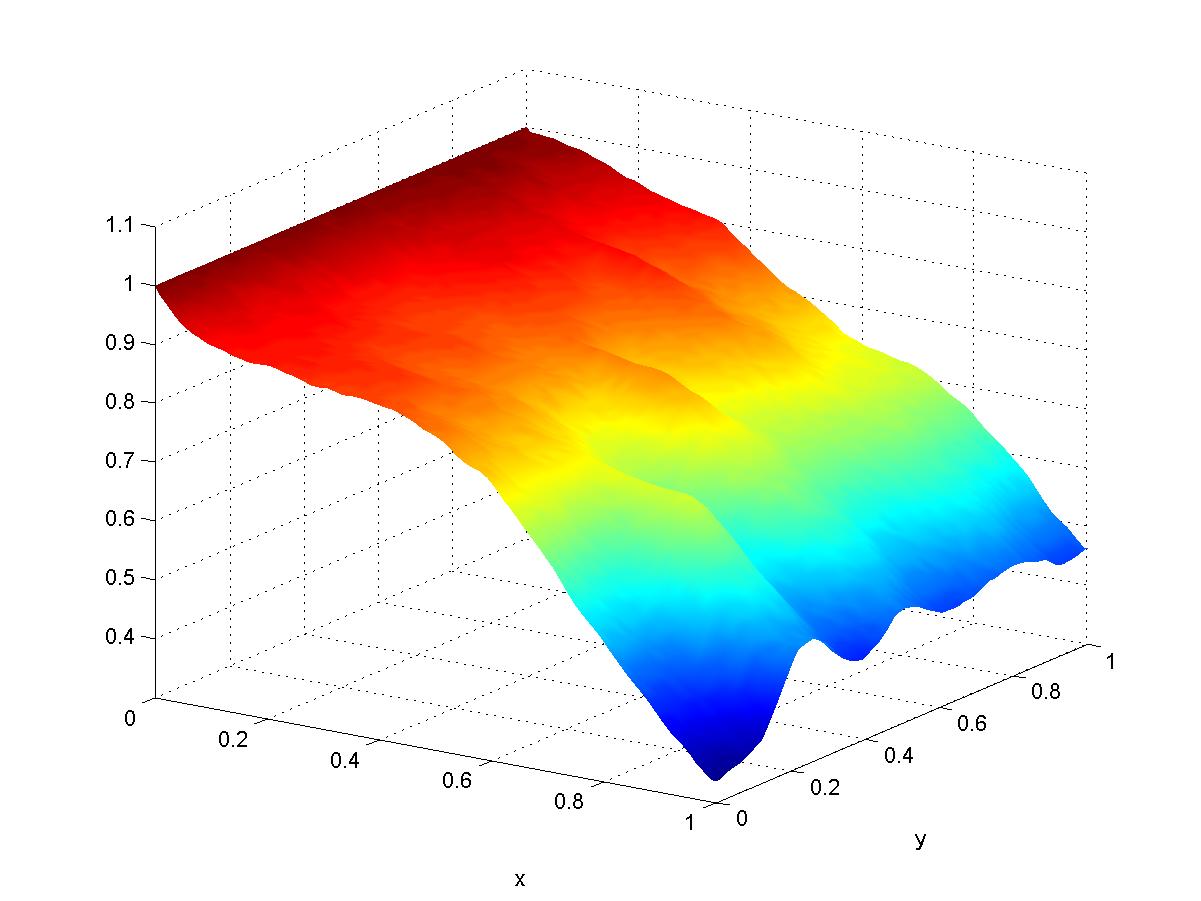}}
 \caption{ (a) Convergence of the root mean square $L^{2}$ norm at $T=1$ as
 a function of $\Dt$ with 1000 realizations with $\Delta x= \Delta y
  = 1/100$, $X_{0}=0$ for \eqref{adr}. The noise is white in time and in $H^\beta$ in space
 $\beta \in\{1,2\}$ (with $\epsilon=0.05$ in \eqref{noise2}). 
 The temporal order of convergence in time is $1/3$. 
 In (b) we plot the streamlines of the velocity field found from solving \eqref{couple1}. 
 In (c) we show a sample reference solution for $\beta=1$ with  $\Dt=1/7680$ while (d) displays 
the mean of 1000 realizations of reference solutions for $\beta=1$. }   
 \label {FIG022}
\end{figure}

Indeed computing the eigenfunctions and eigenvalues of $A$
with  this  mixed Neumann-Dirichlet boundary conditions is expensive. Let's examine the boundary condition and put the problem in an equivalent abstract setting as \eqref{adr}.
Using the trace operator (see \cite{EP}) in Green's theorem yields
\begin{eqnarray}
  dX=(AX +F_{1}(X)+ b(X))dt + d W,
\end{eqnarray}
where for $v\in  H^{1}(\Omega)$ 
$$(A u, v)= -\int_{\Omega} D \nabla u \nabla v\, dx,$$
and 
$$(bu,v)=\int_{\Gamma}\dfrac {\partial  u}{\partial\nu_{A} }
\gamma_{0} v \,d \sigma,\qquad \gamma_{0} v = v\mid_{\partial \Omega},\, v\in  H^{1}(\Omega).$$
In this abstract setting, the  linear operator is $A=D\varDelta$  but  with homogeneous Neumann boundary. The explicit expression of $b$ is unknown. 
To deal with high  P\'{e}clet flows we discretize in space using
finite volumes (viewed as the finite element method (see \cite{AtThesis,FV}).
The finite volume method uses finite difference approximation of
$b=\dfrac {\partial }{\partial\nu_{A}}\mid_{\Gamma}$ (see \cite{FV,EP}). 
The nonlinear term is then $F=F_{1}+b$ where
\begin{eqnarray}
F_{1}(u)=-  \nabla \cdot (\mathbf{q} u)-\frac{u}{(u^2+1)}, \quad u \in H^{1}(\Omega).
\end{eqnarray}
We use a heterogeneous medium with three parallel high permeability streaks, 100 times 
higher compared to the other part of the medium. This could represent for example a highly 
idealized fracture pattern. We obtain the Darcy velocity field $\mathbf{q}$  by  solving \eqref{couple1}
with  Dirichlet boundary conditions on 
$\Gamma_{D}^{1}=\left\lbrace 0 ,1 \right\rbrace \times \left[
  0,1\right] $ and Neumann boundary conditions on
$\Gamma_{N}^{1}=\left( 0,1\right)\times\left\lbrace 0,1\right\rbrace $ such that 
\begin{eqnarray*}
p&=&\left\lbrace \begin{array}{l}
1 \quad \text{in}\quad \left\lbrace 0 \right\rbrace \times\left[ 0,1\right]\\
0 \quad \text{in}\quad \left\lbrace 1 \right\rbrace \times\left[ 0,1\right]
\end{array}\right. 
\end{eqnarray*}
and $- k \,\nabla p (\mathbf{x},t)\,\cdot \mathbf{n} =0$ in  $\Gamma_{N}^{1}$.
Provided $\mathbf{q}$ is bounded,
since $b$ is linear, $F_{1}$ and thus $F$ satisfies \assref{assumption4} (b).
We can write the semi-discrete finite volume method as 
\begin{eqnarray}
dX^{h}=(A_{h}X^{h}+P_{h}F_{1}(X^{h})+P_{h}b(X^{h})) +P_{h}P_N dW,
\end{eqnarray}
where here $A_{h}$ is the space discretization of $D \varDelta $ using
only  homogeneous Neumann boundary conditions and $P_{h}b(X^{h})$
comes from the approximation of diffusion flux on the Dirichlet
boundary condition side (see \cite{FV,AtThesis}).
Thus we can form the noise using eigenfunctions of $\varDelta$ 
with full Neumann boundary conditions for the system with mixed boundary conditions.
We use the noise given by \eqref{noise2} and compute the reference solutions using a time step of $\Delta t = 1/7680$. 
Figure \ref{FIG022a} shows the convergence of the modified 
method and the standard semi-implicit method with noise that is
$H^\beta$ in space with $\beta\in\{1,2\}$. 
We observe that the temporal convergence order is close to $1/3$ for all
the schemes. We observe a reduction order of convergence  in
\thmref{th25}. This reduction order is high, this is probably due to
the velocity $\mathbf{q}$.  
Note that we still obtain an improvement in the accuracy over the standard semi-implicit method.
Figure \ref{FIG022b} shows the streamlines of the velocity field, computed from the elliptic 
equation \eqref{couple1}. \figref{FIG022d} shows a sample reference solution with $\beta=1$ and 
Figure \ref{FIG022d} shows  the mean of 1000 realizations  of the
reference solutions also for $\beta=1$.
%

\section{Proofs of the main results} 
\label{sec:proofs}
\subsection{Some preparatory results}
We introduce  the Riesz representation operator  $R_{h}: V \rightarrow V_{h}$ defined by 
\begin{eqnarray}
\label{riesz}
(-A R_{h}v,\chi)=(-A v,\chi)=a(v,\chi)\qquad \qquad v \in V,\; \forall \chi \in V_{h}.
\end{eqnarray}
Under the regularity assumptions on the triangulation and in view of the $V-$ellipticity \eqref{ellip}, 
it is well known (see \cite{lions}) that the following error bound holds for $v\in V\cap H^{r}(\Omega),$
$ r\in\{1,2\}$
\begin{eqnarray}
\label{regulard}
 \Vert R_{h}v-v\Vert +h \Vert
 R_{h}v-v\Vert_{H^{1}(\Omega)} \leq  C h^{r} \Vert
 v\Vert_{H^{r}(\Omega)}.
\end{eqnarray}
It follows  that for $v\in V\cap H^{r}(\Omega)$
$ r\in\{1,2\}$
\begin{equation}
 \Vert R_{h}v -v \Vert  \leq C h^{r}\Vert v \Vert_{H^{r}(\Omega)}.
\end{equation}
Since $P_h$ is the $L^2$ orthogonal projection  and  $ R_{h}v \in V_h$, we therefore have

\begin{equation}
 \Vert P_{h}v -v \Vert = \underset{\chi  \in V_h }{\inf} \Vert \chi-v\Vert  \leq C h^{r}\Vert v \Vert_{H^{r}(\Omega)}.
\label{eq:feError}
\end{equation}
Since  
\begin{equation}
 \Vert P_{h}v -v \Vert \leq  C\Vert v\Vert, \,   v\in  H,
\end{equation}
We therefore have by interpolation theory

\begin{equation}
 \Vert P_{h}v -v \Vert \leq  C h^r \Vert v\Vert_r, \,   v\in  \mathcal{D}((-A)^r),\,\, 0\leq r\leq 2.
 \label{ph}
\end{equation}
This inequality plays a key role in our  convergence proofs. A similar
inequality for the interpolation operator $I_h$ is given in
\cite[Theorem 3.25, Theorem 3.29]{EP} or in \cite[(2.11),
pp.799]{lions}. 
We start by examining the deterministic linear problem.
Find $u \in V$ such that 
\begin{eqnarray}
\label{homog4}
u'=Au \qquad 
\text{given} \quad u(0)=v, \quad t \in (0,T].
\end{eqnarray}
The corresponding semi--discretization in space is : find $u_{h} \in
V_{h}$ such that  $u_{h}'=A_{h}u_{h}$ where $u_{h}^{0}=P_{h}v$.
The  full discretization of (\ref{homog4}) using implicit Euler in time is given by
\begin{equation}
\label{eq:detEuler}
u_{h}^{n+1} = (\text{I}-\Delta t\;A_{h})^{-(n+1)}\;P_{h}v = S_{h,\Delta t}^{n+1}\;P_{h}v.
\end{equation}
We consider the error at $t_n=n\Dt$ and define the operator $T_n$ below
\begin{equation}
u(t_n)-u_h^n = (S(t_n)-(I-\Delta t\,A_{h})^{-n}\;P_{h})v=: T_n v.
\label{eq:form1}
\end{equation}
\begin{lemma}
\label{lemma1}
Let $0\leq \nu \leq r \leq 2$. For  $v \in \mathcal{D}((-A)^{\nu/2})$,
there exists a constant $C>0$ such that the following estimate holds on the numerical approximation to \eqref{homog4} by \eqref{eq:detEuler}
\begin{eqnarray}
\label{form24}
\Vert u(t_{n})-u_{h}^{n}\Vert &=&\Vert T_{n} v\Vert \leq C
 t_{n}^{-\frac{r-\nu}{2}}( h^{r}+\Delta t^{r/2})\Vert v \Vert_{\nu} .
\end{eqnarray}
 \end{lemma}
 \begin{proof}
 The proof can be found in \cite[Lemma 4.3]{kruse}.
 \end{proof}
 
\begin{lemma}\label{lemma2}
Let $X$ be the  mild solution  of (\ref{adr}) given  in (\ref{eq1}), such  that  \eqref{noisenew} of \assref{assumption4} is satisfied for $0<\beta< 2$.
Let  $t_{1}, t_{2} \in [0,T]$, $t_{1}< t_{2}$. Assume that $X_{0} \in
L_{2}(\mathbb{D},\mathcal{D}((-A)^{\beta/2}$ then we have the following.
%
%

(i) If $F$ satisfies the linear growth condition
 $   \Vert F(X)\Vert \leq C \left(1+\Vert X\Vert \right)$,
for $\sigma=\min(
\beta/2,1/2)$ then
\begin{eqnarray*}
  \mathbf{E}\Vert X(t_{2})- X(t_{1}) \Vert^{2} &\leq&  C
  (t_{2}-t_{1})^{2\sigma} \left(\mathbf{E} \Vert
  X_{0}\Vert_{\beta}^{2}+\underset{0\leq s\leq
    T}{\sup} \mathbf{E}\Vert X(s)\Vert^{2}
  +1 \right).
\end{eqnarray*}
(ii)  
If $X$ is a $H^1(\Omega)$-valued process and $F$ satisfies
the linear growth condition
\begin{eqnarray}
\label{linear2}
\Vert F(X) \Vert \leq C \left(1+\Vert X\Vert_{H^1(\Omega)} \right),
\end{eqnarray}
then
\begin{eqnarray*}
\mathbf{E}\Vert X(t_{2})- X(t_{1}) \Vert^{2} &\leq&  C
(t_{2}-t_{1})^{\beta} \left(\mathbf{E} \Vert
 X_{0}\Vert_{\beta}^{2}+\underset{0\leq s\leq
    T}{\sup} \mathbf{E}\Vert X(s)\Vert_{H^1(\Omega)}^{2}
 +1 \right).
\end{eqnarray*}
%
\end{lemma}
 \begin{proof}
 
See \cite[(2.13) of Theorem 2.4]{xia} for the proof of part i).  This proof can easily be updated for part ii) as  we can  bound $\Vert F(X(s)) \Vert$ by \eqref{linear2}. 

\end{proof}
\subsection{Proof of \thmref{th35}} 
\label{sec:th1}
Recall  that
\begin{eqnarray*}
 X(t_{m})&=&S(t_{m})X_{0}+\underset{k=0}{\sum^{m-1}}\int_{t_{k}}^{t_{k+1}}S(t_{m}-s)F(X(s))ds+\int_{0}^{t_{m}} S(t_{m}-s)d W(s)\\
          &=& \overline{X}(t_{m})+ O(t_{m}).
\end{eqnarray*}
We now estimate 
$\left(\mathbf{E}\Vert X(t_{m})-X_{m}^{h}\Vert^{2}\right)^{1/2}$. 
By construction of the approximation from
\eqref{eq:XbarO} and \eqref{eq:Zmh} we have that 
\begin{eqnarray}
\lefteqn{  X(t_{m})-X_{m}^{h}} \nonumber \\
&=&\overline{X}(t_{m}) + O(t_{m})-\left(Z_{m}^{h}+P_{h}P_{N} O(t_{m})\right)\nonumber \\ 
 &=& \left(\overline{X}(t_{m}) -Z_{m}^{h}\right)+
 \left(P_{N}(O(t_{m}))-P_{h}P_{N}(O(t_{m}))\right)+
\left(O(t_{m})-P_{N}(O(t_{m}))\right)\nonumber \\
&=& I +II +III, \label{eq:IIIIII}
\end{eqnarray}
and $Z_{m}^{h}$ is given by \eqref{eq:Zmh}. 
Then
$$\left(\mathbf{E}\Vert X(t_{m})-X_{m}^{h}\Vert^{2}\right)^{1/2}
\leq \left(\mathbf{E}\Vert I\Vert^{2}\right)^{1/2}
+\left(\mathbf{E}\Vert II\Vert^{2}\right)^{1/2} +\left(\textbf{E}\Vert
 III\Vert^{2}\right)^{1/2} $$
and we estimate each term. Since the first term will require the most work
we estimate $II$ and $III$ first.
Let us examine $\left(\mathbf{E}\Vert II\Vert^{2}\right)^{1/2}$. To
do this we use the finite element estimate \eqref{ph},
the regularity of the noise
and the fact that $P_{N}$ is bounded. 
Then for $0\leq \beta \leq 2 $ , if $O(t_{m}) \in \mathcal{D}((-A)^{\frac{\beta}{2}})$   we have   
$$  \mathbf{E}\Vert II\Vert^{2}  
\leq  C h^{2\beta} \mathbf{E} \Vert P_{N}(O(t_{m}))\Vert_{\beta}^2 
\leq C h^{2\beta} \mathbf{E} \Vert O(t_{m})\Vert_{\beta}^2.$$
Using the Ito isometry, \cite[Lemma 3.2]{xia} and  \eqref{noisenew} of \assref{assumption4} yields
\begin{eqnarray*}
\mathbf{E}\Vert II\Vert^{2} &\leq& 
C h^{2 \beta} \int_{0}^{t_{m}} \Vert (-A)^{\beta/2}S(t_{m}-s) Q^{1/2}\Vert_{\mathcal{L}_2(H)}^{2}ds \\
&\leq & C h^{2 \beta} \Vert (-A)^{(\beta-1)/2}Q^{1/2}\Vert_{\mathcal{L}_2(H)}^2.
\end{eqnarray*}
For the third term $III$ we have 
\begin{eqnarray*}
\mathbf{E}\Vert III\Vert^{2} &=& 
\mathbf{E} \Vert (\text{I}-P_{N})O(t_{m}) \Vert^{2} \\
&=& 
\mathbf{E} \Vert (\text{I}-P_{N})(-A)^{-\beta/2} (-A)^{\beta/2}O(t_{m})\Vert^{2},
\end{eqnarray*}
and so using  again \cite[Lemma 3.2]{xia},  \eqref{noisenew} of \assref{assumption4} 
\begin{eqnarray*}
\mathbf{E}\Vert III\Vert^{2}  &\leq&  
\Vert (\text{I}-P_{N})(-A)^{-\beta/2} \Vert_{L(L^{2}(\Omega))}^2  \mathbf{E} \Vert
(-A)^{\beta/2} O(t_{m})\Vert^2\\
& \leq&  \Vert (-A)^{(\beta-1)/2}Q^{1/2}\Vert_{\mathcal{L}_2(H)}^2 \left( \underset { j \in
    \mathbb{N}^{d} \backslash \mathcal{I}_{N}} {\inf} \lambda_{j}\right)^{-\beta}\\
 & \leq&   C \left( \underset { j \in
    \mathbb{N}^{d} \backslash \mathcal{I}_{N}} {\inf} \lambda_{j}\right)^{-\beta}.
\end{eqnarray*}
To estimate $I$, we follow the approach in \cite{xia} 
\begin{eqnarray*}
 I&=&S(t_{m})X_{0}-S_{h,\Delta t}^{(m)} P_{h}X_{0}\\
 &+& \underset{k=0}{\sum^{m-1}}\int_{t_{k}}^{t_{k+1}}S(t_{m}-s)F(X(s))ds-\underset{k=0}{\sum^{m-1}}\left(\int_{t_{k}}^{t_{k+1}} S_{h,\Delta t}^{(m-k)}P_{h}F(X_{k}^{h})ds\right)\\
&=&I_1+I_2,
\end{eqnarray*}
then
$$\mathbf{E} \Vert I\Vert^{2} \leq 2 \left( \mathbf{E} \Vert
  I_1\Vert^{2}+ \mathbf{E} \Vert I_2\Vert^{2}\right).$$
According to \lemref{lemma1},
for $X_{0}\in L_{2}\left(\mathbb{D}, \mathcal{D}((-A)^{\beta/2})\right),\, 0< \beta < 2 $  we have
 \begin{eqnarray}
 \left(\mathbf{E} \Vert I_1 \Vert^{2}\right)^{1/2} \leq C ( h^{\beta}+\Delta t^{\beta/2}). 
\end{eqnarray}
To estimate $\mathbf{E} \Vert I_2\Vert^2$  we have 
\begin{eqnarray*}
 I_2&=&\underset{k=0}{\sum^{m-1}}\int_{t_{k}}^{t_{k+1}}\left(S(t_{m}-t_k)-S_{h,\Delta t}^{(m-k)} P_h\right) F(X(s))ds \\
  & &+ \underset{k=0}{\sum^{m-1}}\int_{t_{k}}^{t_{k+1}}\left(S(t_{m}-s)-S(t_m-t_k)\right) F(X(s))ds \\
 & & +\underset{k=0}{\sum^{m-1}}\int_{t_{k}}^{t_{k+1}} S_{h,\Delta t}^{(m-k)} P_h \left(F(X(s)-F(X(t_k))\right)ds \\
 & & + \underset{k=0}{\sum^{m-1}}\int_{t_{k}}^{t_{k+1}} S_{h,\Delta t}^{(m-k)} P_{h}\left(F(X(t_k))-F(X_{k}^{h})\right)ds\\
 &=& I_2^{1}+ I_2^{2}+ I_2^{3}+I_2^{4}.
\end{eqnarray*}
Using again \lemref{lemma1},  for $0<\beta<2$ we have
\begin{eqnarray*}
\lefteqn{\left(\mathbf{E} \Vert I_{2}^1\Vert^2\right)^{1/2} }\\
&\leq& \underset{k=0}{\sum^{m-1}}\int_{t_{k}}^{t_{k+1}}\left( \mathbf{E} \Vert \left(S(t_{m}-t_k)-S_{h,\Delta t}^{(m-k)} P_h\right) F(X(s))\Vert^2\right)^{1/2}ds\\
&\leq&  C (h^{(\beta)}+\Delta t^{\beta/2} \left(\underset{k=0}{\sum^{m-1}}\int_{t_{k}}^{t_{k+1}} (t_m-t_k)^{-\beta/2} ds \right)^{1/2}\left(1 + \underset{0\leq s \leq T}{\sup} \mathbf{E} \Vert (X(s))\Vert^{2}\right)^{1/2}\\
                                &\leq&  C (h^{\beta}+\Delta t^{\beta/2}) \left(\int_{0}^{t_{m}} (t_m-t_k)^{-\beta/2} ds\right)^{1/2}\\
                                &\leq& C(h^{\beta}+\Delta t^{\beta/2}).
\end{eqnarray*}
Let us estimate $\left(\mathbf{E} \Vert I_2^{2}\Vert^{2}\right)^{1/2}$. 
By \propref{prop1}, for $0\leq t_{1}< t_{2}\leq T$ we have 
\begin{eqnarray*}
\label{ineq}
\Vert S(t_{2})-S(t_{1})\Vert_{L(L^{2}(\Omega))}&= &\Vert (-A)S(t_{1})
(-A)^{-1}\left(\mathbf{I}-S(t_{2}-t_{1})\right)\Vert_{L(L^{2}(\Omega))} \nonumber\\
& \leq& \dfrac{(t_{2}-t_{1})}{t_{1}}.
\end{eqnarray*}
Splitting the estimation and using \eqref{ineq} in the second part yields
\begin{eqnarray*}
\left(\mathbf{E} \Vert I_2^{2} \Vert^{2}\right)^{1/2} &\leq&\left(
\underset{k=0}{\sum^{m-1}}\int_{t_{k}}^{t_{k+1}}\Vert S( t_{m}-s
)-S(t_{m} -t_{k} ) \Vert_{L(L^{2}(\Omega))} ds \right)\\
&&  \times\left( \underset{0\leq s\leq
    T}{\sup}\mathbf{E} \Vert F(X(s))\Vert^{2}\right)^{1/2}\\ 
&\leq& C \left( \int_{t_{m-1}}^{t_{m}}\Vert S( t_{m}-s
)-S(t_{m} -t_{m-1} ) \Vert_{L(L^{2}(\Omega))} ds \right.\\
&& \left. +
\underset{k=0}{\sum^{m-2}}\int_{t_{k}}^{t_{k+1}}
\left(\dfrac{s-t_{k}}{t_{m}-s}\right)ds\right).
\end{eqnarray*}
Using the fact that the operator $S(t)$ is bounded yields
\begin{eqnarray*}
  \int_{t_{m-1}}^{t_{m}}\Vert S( t_{m}-s
)-S(t_{m} -t_{m-1} ) \Vert_{L(L^{2}(\Omega))} ds \leq  C \Delta t,
\end{eqnarray*}
and then
\begin{eqnarray*}
\left(\mathbf{E} \Vert I_2^{2} \Vert^{2}\right)^{1/2} &\leq&
C \left(\Delta t+  \underset{k=0}{\sum^{m-2}} \left( (m-k-1)\Delta
t\right)^{-1}\int_{t_{k}}^{t_{k+1}}\left(s-t_{k}\right) ds \right)\\ 
&\leq& C \left(\Delta t+  \Delta t \underset{k=0}{\sum^{m-2}} \left(
 m-k-1\right)^{-1}\right). 
 \end{eqnarray*}
Noting that the sum above is bounded by $\ln(M)$ we have 
\begin{eqnarray*}
 \left(\textbf{E} \Vert I_2^{2} \Vert^{2}\right)^{1/2} &\leq&  C (\Delta t+  \Delta t | \ln(\Delta t) |).
\end{eqnarray*}
The estimation of $\mathbf{E} \Vert I_{2}^3\Vert^2$  follows the
one in \cite{xia} using \lemref{lemma2}. 
We provide some keys steps, the main difference comes from the
different space discretization. By Taylor expansion 
\begin{eqnarray*}
 \lefteqn {(X(s)-F(X(t_k))}\\
 &=& F'(X(t_{k}))(X(s)-X(t_k))+RF\\
 &=&  F'(X(t_{k})) \left((S(s-t_k)-I)X(t_k)+  \int_{t_k}^s (S(s-\tau) F(X(\tau))d\tau+ \int_{t_k}^s (S(s-\tau)dW(\tau)\right) \nonumber
\end{eqnarray*}
where $RF= \int_0^1 F''(X(t_{k}))+r(X(s)-X(t_{k}))(
X(s)-X(t_{k}),X(s)-X(t_{k}))dr$. Then,
\begin{eqnarray*}
  \lefteqn{\left(\mathbf{E} \Vert I_{2}^3\Vert^2\right)^{1/2}}\\
&=&\left(\mathbf{E} \Big \Vert \underset{k=0}{\sum^{m-1}}\int_{t_{k}}^{t_{k+1}} S_{h,\Delta t}^{(m-k)} P_h F'(X(t_{k})) \left((S(s-t_k)-I)X(t_k)\right)ds\Big \Vert^2\right)^{1/2}\\
&+& \left(\mathbf{E} \Big \Vert \underset{k=0}{\sum^{m-1}}\int_{t_{k}}^{t_{k+1}} S_{h,\Delta t}^{(m-k)} P_h F'(X(t_{k}))\int_{t_k}^s (S(s-\tau) F(X(\tau))d\tau ds\Big \Vert^2\right)^{1/2}\\
&+& \left(\mathbf{E} \Big \Vert \underset{k=0}{\sum^{m-1}}\int_{t_{k}}^{t_{k+1}} S_{h,\Delta t}^{(m-k)} P_h F'(X(t_{k}))\int_{t_k}^s (S(s-\tau)dW(\tau)ds\Big \Vert^2\right)^{1/2}\\
&+& \left(\mathbf{E} \Big \Vert \underset{k=0}{\sum^{m-1}}\int_{t_{k}}^{t_{k+1}} S_{h,\Delta t}^{(m-k)}P_h RF ds\Big \Vert^2\right)^{1/2}\\
&=& I_{2 1}^3+I_{2 2}^3+I_{2 3}^3+I_{2 4}^3.
\end{eqnarray*}
 Let us estimate $ I_{2 1}^3$. \assref{assumption4} yields
 \begin{eqnarray*}
 \lefteqn{I_{2 1}^3} \\
 &\leq&  \underset{k=0}{\sum^{m-1}}\int_{t_{k}}^{t_{k+1}}\left(\mathbf{E} \Big\Vert S_{h,\Delta t}^{(m-k)} P_h F'(X(t_{k})) \left((S(s-t_k)-I)X(t_k)\right)\Big \Vert^2\right)^{1/2}ds\\
   &\leq&  \underset{k=0}{\sum^{m-1}}\int_{t_{k}}^{t_{k+1}} \Vert
          S_{h,\Delta t}^{(m-k)}
          (-A_h)^{\delta/2}\Vert_{L(L^2(\Omega))} \\
& & \times \left(\mathbf{E} \Big \Vert (-A_h)^{-\delta/2}P_h F'(X(t_{k})) \left((S(s-t_k)-I)X(t_k)\right)\Big \Vert^2\right)^{1/2}ds
  \end{eqnarray*}
The operator $S_{h,\Delta t}$ satisfies the smoothing properties analogous  to
$S(t)$  independently of $h$ (see for example \cite{Stig,ElliottLarsson}),
we find for $t_{m}=m\Delta t >0$
\begin{eqnarray*}
 \label{smooth}
  S_{h,\Delta t}^{m} (-A_h)^{\delta/2} \leq C t_m^{-\delta/2}.
\end{eqnarray*}
Note that for $u\in H=L^2(\Omega)$
\begin{eqnarray}
 \Vert (-A_h)^{-\frac{\delta}{2}} P_h u\Vert=\underset{v_h \in V_h}{\sup} \dfrac{((-A_h)^{-\delta/2} P_h u,v_h)}{\Vert v_h \Vert}= \underset{v_h \in V_h}{\sup} \dfrac{(P_h u, (-A_h^*)^{-\delta/2} v_h )}{\Vert v_h \Vert}.
\end{eqnarray}
Using the definition of $P_h$ we have 
\begin{eqnarray}
\label{keyextimation}
 \Vert (-A_h)^{-\frac{\delta}{2}} P_h u\Vert&=& \underset{v_h \in
                                                V_h}{\sup} \dfrac{(u,
                                                (-A_h^*)^{-\delta/2}
                                                v_h )}{\Vert v_h
                                                \Vert}=\underset{w_h
                                                \in V_h}{\sup}
                                                \dfrac{(u, w_h)}{\Vert
                                                (-A_h^*)^{\delta/2}
                                                v_h \Vert} \nonumber \\ 
 &\leq& C
 \underset{w_h \in V_h}{\sup} \dfrac{(u, w_h )}{\Vert v_h \Vert_{\delta}}= \Vert u\Vert_{-\delta} = \Vert (-A)^{-\delta/2} u\Vert.
\end{eqnarray}
%
Using \eqref{keyextimation} and \eqref{smooth} yields
\begin{eqnarray*}
I_{2 1}^3&\leq& 
   C \underset{k=0}{\sum^{m-1}}\int_{t_{k}}^{t_{k+1}}(t_m-t_k)^{-\delta/2}\left(\mathbf{E} \Big \Vert(-A)^{-\delta/2}F'(X(t_{k})) \left((S(s-t_k)-I)X(t_k)\right)\Big \Vert^2\right)^{1/2}ds
 \end{eqnarray*}
 The rest of the estimation follows \cite[$I_{3 1}$]{xia} using  \lemref{lemma2}, and we have
 $$ I_{2 1}^3 \leq C(\Delta t^{\min(\beta,1)}).$$
 It is obvious that  $ I_{2 2}^3 \leq C \Delta t.$ The estimation of
 $I_{2 3}^3$ follows from \cite[$I_{3 3}$]{xia} by 
 replacing $P_N$ by $P_h$ and $E(t_m-s)$ by $ S_{h,\Delta t}^{(m-k)}$
 to get $ I_{2 3}^3 \leq C \Delta t^{\min(1+\beta,2)}.$
 The estimation of  $I_{2 4}^3$ follows from \cite[$I_{3 4}$]{xia} by
 replacing $P_N$ by $P_h$ and $E(t_m-s)$ by $ S_{h,\Delta t}^{(m-k)}$
 and we finally have 
 $I_{2 4}^3 \leq  \Delta t^{\min(\beta,1)}.$\\
 Putting the estimates on $I_{21}^3$, $I_{22}^3$, $I_{23}^3$ and $I_{24}^3$ together we see that 
\begin{eqnarray}
\left(\ \mathbf{E} \Vert I_2^2 \Vert^{2}\right)^{1/2} \leq C \Delta t^{\min(\beta,1)}. 
\end{eqnarray}
For $\mathbf{E} \Vert I_{2}^3\Vert^2$, we obviously have
\begin{eqnarray}
 \left(\mathbf{E} \Vert I_{2}^3\Vert^2\right)^{1/2} \leq  C \underset{k=0}{\sum^{m-1}}\int_{t_{k}}^{t_{k+1}}\left(\mathbf{E} \Vert X(t_k)- X_{k}^{h}\Vert^2\right)^{1/2}ds.
\end{eqnarray}
Combining the estimates of $ \mathbf{E} \Vert I\Vert^{2}$ and $
\mathbf{E} \Vert II\Vert^{2}$ and applying the discrete Gronwall
lemma completes the proof.
\subsection{Proof of \thmref{th25}}
\label{sec:th2}
The estimation of $I $ and $II$ is the same as before.
We now estimate the term $I$ from \eqref{eq:IIIIII} when there is
non-zero advection.
As above we have 
\begin{eqnarray}
I &=&  T_{m}X_{0} + \underset{k=0}{\sum^{m-1}} \int_{t_{k}
 }^{t_{k+1}} S( t_{m}-s) F(X(s))-S_{h,\Delta t}^{(m-k)}
P_{h}F(Z_{k}^{h}+P_{h}P_{N}O(t_{k}))ds \nonumber \\  
 &=& T_{m}X_{0} + \underset{k=0}{\sum^{m-1}} \int_{t_{k}}^{t_{k+1}} S_{h,\Delta t}^{(m-k)}P_{h}(F(X(t_{k}))
-F(Z_{k}^{h}+P_{h}P_{N}O(t_{k}))))ds\nonumber \\ 
&& +\underset{k=0}{\sum^{m-1}} \int_{t_{k}}^{t_{k+1}}
S_{h,\Delta t}^{(m-k)}P_{h}(F(X(s))-F (X(t_{k})))ds\nonumber \\ 
&& +\underset{k=0}{\sum^{m-1}} \int_{t_{k}}^{t_{k+1}}(S( t_{m}-
t_{k} )-S_{h,\Delta t}^{(m-k)}P_{h})F(X(s))ds\nonumber \\ 
&& +\underset{k=0}{\sum^{m-1}} \int_{t_{k}}^{t_{k+1}} (S(t_{m} 
-s )-S(t_{m} -t_{k}))F(X(s))ds\nonumber \\ 
&=& I_{1}+I_{2}+I_{3}+I_{4}+I_{5}.
\label{eq:I1toI5}
 \end{eqnarray}
The estimation of $I_1$ is the same as for \thmref{th35}.  
If $F$ satisfies \assref{assumption4}(b), 
using \eqref{keyextimation} with  $\delta=1$ yields
\begin{eqnarray*}
\lefteqn{(\mathbf{E} \Vert I_{2} \Vert^2)^{1/2}}\\ &\leq &
\underset{k=0}{\sum^{m-1}}\int_{t_{k}}^{t_{k+1}} \left(\mathbf{E} \Vert
S_{h,\Delta t}^{(m-k)} P_{h}\left(F(X(t_{k}))-F(Z_{k}^{h}+P_{h}P_{N}O(t_{k}))\right)\Vert^{2}\right)^{1/2}ds\\ 
&=& \underset{k=0}{\sum^{m-1}}\int_{t_{k}}^{t_{k+1}} \left(\mathbf{E} \Vert
S_{h,\Delta t}^{(m-k)} (-A_h)^{1/2} (-A_h)^{-1/2} P_{h}\left(F(X(t_{k}))-F(Z_{k}^{h}+P_{h}P_{N}O(t_{k}))\right)\Vert^{2}\right)^{1/2}ds\\ 
&\leq &  C \underset{k=0}{\sum^{m-1}}\int_{t_{k}}^{t_{k+1}} (t_{m}-t_k)^{-1/2}\left(\mathbf{E} \Vert F(X(t_{k}))-F(Z_{k}^{h}+P_{h}P_{N}O(t_{k}))\Vert_{-1}^{2}\right)^{1/2}ds\\
 &\leq & C \underset{k=0}{\sum^{m-1}}\int_{t_{k}}^{t_{k+1}}(t_{m}-t_k)^{-1/2}\left(\mathbf{E}\Vert X(t_{k})-X_{k}^{h}\Vert^{2}\right)^{1/2}ds.
\end{eqnarray*}

Let us estimate $(\mathbf{E} \Vert I_{3} \Vert^{2})^{1/2}$. Once again using the Lipschitz  condition and smoothing
property of $ S_{h, \Delta t}$
\begin{align*}
(\mathbf{E} \Vert I_{3} \Vert^{2})^{1/2}
&\leq  C \underset{k=0}{\sum^{m-1}}\int_{t_{k}}^{t_{k+1}}
(t_{m}-t_k)^{-1/2}\left(\mathbf{E}  \Vert F(X(s))-F(X(t_{k}))\Vert_{-1}\right)^{1/2}ds\\ 
&\leq  C\underset{k=0}{\sum^{m-1}}\int_{t_{k}}^{t_{k+1}}
(t_{m}-t_k)^{-1/2}(\mathbf{E}  \Vert
X(s)-X(t_{k})\Vert^{2})^{1/2}ds
\end{align*}
Since 
$$\underset{k=0}{\sum^{m-1}}\int_{t_{k}}^{t_{k+1}}(t_{m}-t_k)^{-1/2}ds\leq 2
\sqrt{ T}.$$
Then if  $X_{0} \in  L_{2}(\mathbb{D}, \mathcal{D}((-A)^{\beta/2})),\, 0<\beta<2$ 
such that $X(t)$ is a $H^1(\Omega)$- valued  process
\begin{eqnarray*}
  (\mathbf{E} \Vert I_{3} \Vert^{2})^{1/2} \leq C
  (\Delta t)^{\min(\beta/2,1/2)}\left(\mathbf{E} \Vert
    X_{0}\Vert_{\beta}^{2}+\underset{0\leq s\leq
        T}{\sup} \mathbf{E} \Vert X(s)\Vert_{H^1(\Omega)}^{2}+1\right)^{1/2},
\end{eqnarray*}
From \lemref{lemma1}, for $0\leq \beta <2$ such  that $X(t)$ is a $H^1(\Omega)$- valued  process,  we have 
\begin{eqnarray*}
 (\mathbf{E} \Vert I_{4} \Vert^{2})^{1/2} &\leq& \underset{k=0}{\sum^{m-1}}\int_{t_{k}}^{t_{k+1}} \left(\mathbf{E} \Vert T_{m}(t_{m}-t_k) F(X(s))\Vert^2\right)^{1/2} ds\\
 & \leq & C (h^{\beta}+\Delta t^{\beta/2}) \left(\underset{k=0}{\sum^{m-1}}\int_{t_{k}}^{t_{k+1}}  (t_{m}-t_k)^{-\beta/2}ds \right) \left(\underset{0\leq s\leq T}{\sup} \mathbf{E} \Vert F(X(s)) \Vert^{2}\right)^{1/2}\\
& \leq & C(h^{\beta}+\Delta t^{\beta/2}) \left(1+\left(\underset{0\leq s\leq T}{\sup} \mathbf{E} \Vert  X(s) \Vert_{H^1(\Omega)}^{2}\right)^{1/2}\right)\\
& \leq & C (h^{\beta}+\Delta t^{\beta/2}) 
\end{eqnarray*}
As in the estimation of $I_2^2$ in the previous section, if  $X(t)$ is a $H^1(\Omega)$- valued  process we have
\begin{eqnarray*}
  (\mathbf{E} \Vert I_{5} \Vert^{2})^{1/2} &\leq&C(\Delta t+ \Delta t \vert \ln (\Delta t)\vert).
\end{eqnarray*}


Combining our estimates $\left(\mathbf{E} \Vert
  I\Vert^{2}\right)^{1/2},\left(\mathbf{E} \Vert
  II\Vert^{2}\right)^{1/2}$ and $\left(\mathbf{E} \Vert
  III\Vert^{2}\right)^{1/2}$ and using the discrete Gronwall lemma 
concludes the proof.

\subsection{Proof of Proposition \ref{Prop:QEVALS}} 
\label{sec:Prop:QEVALS}
Let  $b$ and $\lambda$ be two real numbers, then the following result holds
\begin{eqnarray}
\label{exp}
 \int_{-\infty}^{+ \infty}\exp \left(-\dfrac{\pi}{4}\left(\frac{ x^{2}
     }{b^{2}}\right) \right) \cos(\lambda x) dx =  2 b\exp\left[
   -\dfrac{1}{\pi}\left(\lambda b\right)^{2} \right]. 
\end{eqnarray}
Recall \cite{DaPZ} that the covariance operator $Q$ may be defined for $f \in
L^{2}(\Omega)$ by  
\begin{eqnarray*}
 Qf(x) = \int_{\Omega}C_{r}(x,y) f(y)dy.
\end{eqnarray*}
Indeed  we have 
\begin{eqnarray*}
\label{eig}
  \lefteqn{4 b_{1} b_{2} \int_{0}^{L_{1}} \int_{0}^{L_{2}}
  C_{r}((x_{1},y_{1});(x_{2},y_{2})) \cos(\lambda_{i}^{(1)} x_{2}
  )\cos(\lambda_{j}^{(2)} y_{2} )dy_{2} dx_{2}}\\ 
&= & \Gamma\;\int_{0}^{L_{1}} \exp
\left(-\dfrac{\pi}{4}\left(\dfrac{\left( x_{2}-x_{1}\right)^{2}
    }{b_{1}^{2}}\right) \right) \cos(\lambda_{i}^{(1)} x_{2} )dx_{2} \\ 
& &\times \int_{0}^{L_{2}} \exp \left(-\dfrac{\pi}{4}\left[\dfrac{\left(
        y_{2}-y_{1}\right)^{2} }{b_{2}^{2}}\right] \right)
\cos(\lambda_{j}^{(2)} y_{2} ) dy_{2}\\
&=& \Gamma\; \int_{-x_{1}}^{L_{1}-x_{1}}\exp
\left(-\dfrac{\pi}{4}\left(\dfrac{ x^{2} }{b_{1}^{2}}\right) \right)
\cos(\lambda_{i}^{(1)} (x+x_{1} ))dx \\ 
& & \times \int_{-y_{1}}^{L_{2}-y_{1}}
\exp \left(-\dfrac{\pi}{4}\left(\dfrac{ x^{2} }{b_{2}^{2}}\right)
\right) \cos(\lambda_{j}^{(2)}(x+y_{1})dx).
\end{eqnarray*}
For $b_{i}\ll L_{i}$, because of the strong decay, we approximate the integral in the finite domain by the integral in 
infinite domain where we can evaluate exactly 
\begin{eqnarray*}
\lefteqn{4 b_{1} b_{2} \int_{0}^{L_{1}} \int_{0}^{L_{2}}
  C_{r}((x_{1},y_{1});(x_{2},y_{2})) \cos(\lambda_{i}^{(1)} x_{2}
  )\cos(\lambda_{j}^{(2)} y_{2} )dy_{2} dx_{2}}\\ 
&\approx& \Gamma \;\int_{-\infty}^{+ \infty}\exp
\left(-\dfrac{\pi}{4}\left(\dfrac{ x^{2} }{b_{1}^{2}}\right) \right)
\cos(\lambda_{i}^{(1)} (x+x_{1} ))dx \\ && \times \int_{-\infty}^{+ \infty}
\exp \left(-\dfrac{\pi}{4}\left(\dfrac{ x^{2} }{b_{2}^{2}}\right)
\right) \cos(\lambda_{j}^{(2)}(x+y_{1}))dx\\ 
&=& 4 b_{1} b_{2} \cos(\lambda_{i}^{(1)}\,x_{1})
\cos(\lambda_{j}^{(2)}\,y_{1})\,\Gamma\,\exp\left(
  -\dfrac{1}{\pi}\left((\lambda_{i}^{(1)}b_{1})^{2}+(\lambda_{j}^{(2)}b_{2})^{2}\right) \right).   
\end{eqnarray*}
It is important to notice that we have used the fact that 
\begin{eqnarray*}
  \int_{-\infty}^{+ \infty}\exp \left(-\dfrac{\pi}{4}\left(\dfrac{
        x^{2} }{b_{i}^{2}}\right) \right) \cos(\lambda_{j}^{(i)}x)dx
  &= & 2 b_{i}\exp\left[ -\dfrac{1}{\pi}\left((\lambda_{j}^{(i)}
      b_{i})^{2}\right) \right] \;\;\;\; i \in \left\lbrace
    1,2\right\rbrace  
\end{eqnarray*}
by \eqref{exp} and
\begin{eqnarray*}
\int_{-\infty}^{+ \infty}\exp \left(-\dfrac{\pi}{4}\left(\dfrac{ x^{2} }{b_{i}^{2}}\right) \right) \sin(\lambda_{j}^{(i)}x)dx &=&0 
\end{eqnarray*}
because the integrand is an odd function.
Then the corresponding values of $\left\lbrace q_{i,j}\right\rbrace_{i+j > 0}$ in the representation (\ref{eq:W1}) is given by
\begin{eqnarray*}
 q_{i,j}= \Gamma \exp\left[ -\dfrac{1}{2 \pi}\left((\lambda_{i}^{(1)}b_{1})^{2}+(\lambda_{j}^{(2)}b_{2})^{2}\right) \right].
\end{eqnarray*}

\subsection*{Acknowledgements}
We would like to thank Dr A. Jentzen for very useful discussions at an
early stage of this paper. These were made possible through an 
{\sc  arc--daad} grant number 1333. 
Antoine  Tambue  was also funded by the Overseas Research Students Awards Scheme (ORSAS), Heriot Watt University and  Robert Bosch Stiftung through the AIMS ARETE chair programme.

\end{document}